\documentclass[onecolumn,final]{elsart3p}
\usepackage[a4paper, total={6in, 8in}]{geometry}

\usepackage{algorithm} 
\usepackage{algpseudocode} 
\usepackage{appendix}
\usepackage{cite}
\usepackage{tabularx}

\newtheorem{theo}{Theorem}[section]

\newtheorem{gass}{Assumption G\hspace{-1.mm}}

\newcounter{exnum}[section]

\newcommand{\Dtriangle}{\calD_{\rm{Triangle}}}
\newcommand{\Dmaze}{\calD_{\rm{Maze}}}
\newcommand{\Dstar}{\calD_{\rm{Star}}}
\newcommand{\Djenga}{\calD_{\rm{Jenga}}}
\newcommand{\Dslices}{\calD_{\rm{Slices}}}
\newcommand{\Dulike}{\calD_{\rm{Ulike}}}
\newcommand{\DjengaM}{\calD_{\rm{Jenga4}}}
\newcommand{\DslicesM}{\calD_{\rm{Slices4}}}
\newcommand{\DulikeM}{\calD_{\rm{Ulike4}}}

\usepackage{amsmath} % AMS equation environments
\usepackage{amssymb} % AMS symbol fonts (e.g., \boldsymbol{})
\usepackage{bbm}
\usepackage{booktabs}
\usepackage{bm}
\usepackage{chapterbib} % Bibtex
\usepackage{color}
\usepackage{comment}
\usepackage{epsfig}
\usepackage{graphics}
\usepackage{graphicx}
\usepackage{multicol}
\usepackage{multirow}
\usepackage{overpic}
\usepackage{pgfplots}
\usepackage{pgfplotstable}
\pgfplotsset{compat=newest}
\usepackage{rotating}
\usepackage{siunitx}
\usepackage{subfigure}
\usepackage{times}
\usepackage{tikz}

\theoremstyle{plain}
\newtheorem{theorem}{Theorem}[section]

\newtheorem{corollary}[theorem]{Corollary}

\theoremstyle{definition}

\theoremstyle{plain}

\def\trait #1 #2 #3 {\vrule width #1pt height #2pt depth #3pt}
\def\fin{\hfill
        \trait .3 5 0
        \trait 5 .3 0
        \kern-5pt
        \trait 5 5 -4.7
        \trait 0.3 5 0
\medskip}

\newcommand{\KER} {\textrm{ker}}
\newcommand{\INTP}{I}%{\footnotesize{\texttt{I}}}

\newcommand{\REAL}{\mathbbm{R}}

\newcommand{\PGRAPH}[1]{\medskip\noindent\textbf{#1}}

\renewcommand{\cv}{\mathbf{c}}

\newcommand{\kv}{\mathbf{k}}

\newcommand{\xv}{\mathbf{x}}
\newcommand{\yv}{\mathbf{y}}

\newcommand{\calD}{\mathcal{D}}

\newcommand{\calH}{\mathcal{H}}
\newcommand{\calI}{\mathcal{I}}

\newcommand{\calM}{\mathcal{M}}
\newcommand{\calN}{\mathcal{N}}

\newcommand{\calT}{\mathcal{T}}

\newcommand{\calV}{\mathcal{V}}

\newcommand{\as}{a}

\newcommand{\cs}{c}

\newcommand{\fs}{f}

\newcommand{\ks}{k}

\newcommand{\qs}{q}

\renewcommand{\ss}{s}

\newcommand{\us}{u}
\newcommand{\vs}{v}
\newcommand{\ws}{w}
\newcommand{\xs}{x}
\newcommand{\ys}{y}

\newcommand{\Cs}{C}

\newcommand{\Fs}{F}

\newcommand{\Ns}{N}

\newcommand{\Ss}{S}

\newcommand{\xsP}{\xs_{\P}}
\newcommand{\ysP}{\ys_{\P}}
\newcommand{\xsE}{\xs_{\E}}
\newcommand{\ysE}{\ys_{\E}}

\newcommand{\PS}[1]{\mathbbm{P}_{#1}}

\newcommand{\HONE}   {H^1}
\newcommand{\HONEzr} {H^1_0}

\newcommand{\LTWO}  {L^2}

\newcommand{\HS}[1] {H^{#1}}

\newcommand{\CS}[1] {C^{#1}}

\renewcommand{\P} {E}% polyhedral element
\newcommand  {\E} {e}% {\textsf{e}}            % edge
\newcommand{\hh}{h}
\newcommand{\Th}{\Omega_{\hh}}

\newcommand{\xvP}{\xv_{\P}}        % element
\newcommand{\xvE}{\xv_{\E}}        % edge
\newcommand{\hP}{\hh_{\P}}

\newcommand{\hE}{\hh_{\E}}

\newcommand{\mP}{\ABS{\P}}

\newcommand{\mE}{\ABS{\E}}

\newcommand{\dV}{\,dV}
\newcommand{\dS}{\,dS}
\newcommand{\dxv}{\,d\xv}

\newcommand{\fsh}{\fs_{\hh}}

\newcommand{\ush}{\us_{\hh}}

\newcommand{\usI}{\us_{\INTP}}

\newcommand{\vsh}{\vs_{\hh}}

\newcommand{\wsh}{\ws_{\hh}}

\newcommand{\wsI}{\ws_{\INTP}}

\newcommand{\Fsh}{\Fs_{\hh}}

\newcommand{\asP}{\as^{\P}}

\newcommand{\ash}{\as_{\hh}}

\newcommand{\ashP}{\as^{\P}_{\hh}}

\newcommand{\nlen}{\hspace{-0.2mm}}

\newcommand{\snorm}  [2]{|#1|_{#2}}

\newcommand{\norm}   [2]{|\nlen|#1|\nlen|_{#2}}

\newcommand{\Tnorm}  [2]{|\nlen|\nlen|#1|\nlen|\nlen|_{{#2}}}

\newcommand{\ABS}    [1]{\left|#1\right|}

\newcommand{\Vhk} {V^{\hh}_{k}}

\newcommand{\Pin}[1]{\Pi^{\nabla}_{#1}}
\newcommand{\Piz}[1]{\Pi^{0}_{#1}}
\newcommand{\PinP}[1]{\Pi^{\nabla,\P}_{#1}}
\newcommand{\PizP}[1]{\Pi^{0,\P}_{#1}}

\newcommand{\SPh} {S^{\P}_{\hh}}
\newcommand{\SPhdd} {S^{\P,\textrm{dd}}_{\hh}}
\newcommand{\SPhtr} {S^{\P,\textrm{tr}}_{\hh}}
\newcommand{\cbot}{c_*}
\newcommand{\ctop}{c^*}

\newcommand{\restrict}[2]{{#1}_{|{#2}}}

\newcommand{\EOD}{\end{document}}

\newcommand{\roundPrecision}{2}
\newcommand{\NDOFS}{N^{\textrm{dofs}}}

\newcommand{\matPin}[1]{\boldsymbol{\Pi}^{\nabla}_{#1}}
\newcommand{\matPiz}[1]{\boldsymbol{\Pi}^{0}_{#1}}

\begin{document}
\begin{frontmatter}
  
  \title{The role of mesh quality and mesh quality indicators in\\
  the Virtual Element Method}
  
  \author[IMATI]{T. Sorgente}
  \author[IMATI]{, S. Biasotti}
  \author[IMATI]{, G. Manzini}
  \author[IMATI]{, M. Spagnuolo}
  
  \address[IMATI]{
    Istituto di Matematica Applicata e Tecnologie Informatiche,
    Consiglio Nazionale delle Ricerche,
    Italy
  }

  % Abstract
  % ----------------------------
  \begin{abstract}
    Since its introduction, the Virtual Element Method (VEM) was shown
    to be able to deal with a large variety of polygons, while
    achieving good convergence rates.
    The regularity assumptions proposed in the VEM literature to guarantee the convergence on a theoretical basis are therefore quite general.
    They have been deduced in analogy to the similar conditions developed in the Finite Element Methods (FEMs) analysis.
    In this work, we experimentally show that the VEM still converges with almost
    optimal rates and low errors in the $\LTWO$ and $\HONE$ norms even
    if we significantly break the regularity assumptions that are used
    in the literature.
    These results suggest that the regularity assumptions proposed so far might be overestimated.
    We also exhibit examples on which the VEM sub-optimally converges
    or diverges.
    Finally, we introduce a mesh quality indicator
    that experimentally correlates the entity of the violation of the
    regularity assumptions and the performance of the VEM solution,
    thus predicting if a dataset is potentially critical for VEM.
  \end{abstract}
  
  
  % \noindent {\bf MSC}: 60J28, 60J85, 92D30, 92D40, 60H10
\end{frontmatter}

\raggedbottom

%% paper
%%\input{sec1_introduction}
%%\input{sec2_vem}
%%\input{sec3_geometrical}
%%\input{sec4_breaking}
%%\input{sec5_metric}
%%\input{sec6_conclusions}
%%\input{app1_error}
%%\input{app2_dataset}

\section{Introduction}
\label{sec:intro}

Finite element methods are very successful in the numerical
treatment of partial differential equations (PDEs), but their
formulation requires an explicit knowledge of the basis functions.
Consequently, they are mostly restricted to meshes with elements
having a simple geometrical shape, such as triangles or quadrilaterals.
This restriction is overcome by polytopal element methods such as the
VEM, which are designed to provide arbitrary
order of accuracy on more generally shaped elements.
In the VEM setting, we partition the computational domain into
polytopal elements and the explicit knowledge of the basis functions
is not required, since the VEM formulation and its practical
implementation is based on suitable polynomial projections that are
always computable from a careful choice of the degrees of freedom.

The VEM was originally formulated
in~\cite{BeiraodaVeiga-Brezzi-Cangiani-Manzini-Marini-Russo:2013} as a
conforming FEM for the Poisson problem by rewriting in a variational
setting the \emph{nodal} mimetic finite difference (MFD)
method~\cite{Brezzi-Buffa-Lipnikov:2009,%
  BeiraodaVeiga-Lipnikov-Manzini:2011,%
  BeiraodaVeiga-Manzini-Putti:2015,%
  Manzini-Lipnikov-Moulton-Shashkov:2017%
} for solving diffusion
problems on unstructured polygonal meshes.
A survey on the MFD method can be found in the review
paper~\cite{Lipnikov-Manzini-Shashkov:2014} and the research
monograph~\cite{BeiraodaVeiga-Lipnikov-Manzini:2014}.
The VEM scheme inherits the flexibility of the MFD method with respect to
the admissible meshes and this feature is well reflected in the many
significant applications that have been developed so far, see, for
example,~\cite{%
  BeiraodaVeiga-Manzini:2014,%
  Benedetto-Berrone-Pieraccini-Scialo:2014,%
  BeiraodaVeiga-Manzini:2015, Berrone-Pieraccini-Scialo-Vicini:2015,%
  Mora-Rivera-Rodriguez:2015,%
  Paulino-Gain:2015,%
  Antonietti-BeiraodaVeiga-Scacchi-Verani:2016,%
  BeiraodaVeiga-Chernov-Mascotto-Russo:2016,%
  BeiraodaVeiga-Brezzi-Marini-Russo:2016a,%
  BeiraodaVeiga-Brezzi-Marini-Russo:2016b,%
  Cangiani-Georgoulis-Pryer-Sutton:2016,%
  Perugia-Pietra-Russo:2016,%
  Wriggers-Rust-Reddy:2016,
  Certik-Gardini-Manzini-Vacca:2018:ApplMath:journal,%%
  Dassi-Mascotto:2018,%
  Benvenuti-Chiozzi-Manzini-Sukumar:2019:CMAME:journal,%
  Antonietti-Manzini-Verani:2019:CAMWA:journal,%%
  Certik-Gardini-Manzini-Mascotto-Vacca:2020%
}.
Because of its origins, the VEM is intimately connected with other
finite element approaches.
The connection between the VEM and finite elements on
polygonal/polyhedral meshes is thoroughly investigated
in~\cite{Manzini-Russo-Sukumar:2014,
  Cangiani-Manzini-Russo-Sukumar:2015, DiPietro-Droniou-Manzini:2018},
between VEM and discontinuous skeletal gradient discretizations
in~\cite{DiPietro-Droniou-Manzini:2018}, and between the VEM and the
BEM-based FEM method
in~\cite{Cangiani-Gyrya-Manzini-Sutton:2017:GBC:chbook}.
The VEM has been extended to convection-reaction-diffusion problems
with variable coefficients
in~\cite{BeiraodaVeiga-Brezzi-Marini-Russo:2016b}.

\medskip
Optimal convergence rates for the virtual element approximations of
the Poisson equation were proved in $\HONE$ and $\LTWO$ norms, see for
instance
\cite{BeiraodaVeiga-Brezzi-Cangiani-Manzini-Marini-Russo:2013,
  Ahmad-Alsaedi-Brezzi-Marini-Russo:2013,CHINOSI2016,
  Beirao-Lovadina-Russo:2016, Brenner:2017:SEV, brenner2018virtual,
  da2020sharper}.
The theoretical results behind the VEM convergence rate involve an estimate of the approximation
error, which is due to both analytical assumptions (interpolation and polynomial projections of the virtual
element functions) and geometrical assumptions (the geometrical shape of
the mesh elements).
There is a general concordance in the literature about the analytical assumptions, but the understanding of which geometrical features of the mesh elements influence the most on the approximation error and the convergence rate, is still an open issue.
Various geometrical (or \textit{regularity}) assumptions have
been proposed to ensure that all elements of any mesh of a given mesh
family in the refinement process are sufficiently regular.
These assumptions guarantee the VEM convergence and optimal estimates
of the approximation error with respect to different norms.
However, as already observed from the very first papers,
cf.~\cite{Ahmad-Alsaedi-Brezzi-Marini-Russo:2013}, the VEM seems to
maintain its optimal convergence rates also when we use mesh families
that do not satisfy the usual geometrical assumptions.

\medskip
As a first contribution of this paper, we overview the geometrical assumptions introduced in the VEM literature to guarantee the convergence.
Then, we define a mesh generation framework that allows us 
to build sequences of meshes (datasets) gradually introducing several pathologies.
The so-generated datasets systematically violate the geometrical
assumptions, and enhance a correlation analysis between 
such assumptions and the VEM performance.
We experimentally show how the VEM presents a good convergence rate on
most examples and only fails in very few situations.
We also provide an indicator of the violation of the geometrical
assumptions, which depends uniquely on the geometry of the mesh elements.
We show a correspondence between this indicator and the performance of the VEM on a given mesh, or mesh family, in terms of approximation
error and convergence rate.
Our work is focused on developing a strategy to evaluate if a given
sequence of meshes is suited to the virtual element
discretization, and possibly to predict the behaviour of the numerical
discretization \emph{before} applying the method.
In this sense, we can consider the approach that we present in this
work as more in an \emph{a priori} than an \emph{a posteriori}
setting.

\medskip
%% OUTLINE
The paper is organized as follows.
In Section~\ref{sec:VEM}, we present the VEM and the convergence
results for the Poisson equation with Dirichlet boundary conditions.
In Section~\ref{sec:geometrical:assumptions}, we detail the geometrical
assumptions on the mesh elements that are used in the literature to
guarantee the convergence of the VEM.
In Section~\ref{sec:breaking}, we present a number of datasets 
which do not satisfy these assumptions, and 
experimentally investigate the convergence of the VEM over them.
In Section~\ref{sec:indicator}, we propose a mesh quality indicator to
predict the behaviour of the VEM over a given dataset.
In Section~\ref{sec:conclusions}, we offer our concluding remarks and
discuss future developments and work.
In Appendices~A and~B, we review the major theoretical results on the
error analysis that are available in the virtual element literature, reporting the geometrical conditions assumed in each result,
and present the algorithmic procedures that we used to build the datasets.

\subsection{Notation and technicalities}
\label{subsec:notation}
We use the standard definition and notation of Sobolev spaces, norms
and seminorms, cf.~\cite{Adams-Fournier:2003}.
Let $k$ be a nonnegative integer number.
The Sobolev space $\HS{k}(\omega)$ consists of all square integrable
functions with all square integrable weak derivatives up to order $k$
that are defined on the open, bounded, connected subset $\omega$ of
$\REAL^{d}$, $d=1,2$.
As usual, if $k=0$, we prefer the notation $\LTWO(\omega)$.
Norm and seminorm in $\HS{k}(\omega)$ are denoted by
$\norm{\cdot}{k,\omega}$ and $\snorm{\cdot}{k,\omega}$, 
while for the inner product in $\LTWO(\omega)$ we prefer 
the integral notation.
We denote the space of polynomials of degree less than or equal to
$k\geq0$ on $\omega$ by $\PS{k}(\omega)$ and conventionally assume
that $\PS{-1}(\omega)=\{0\}$.
In our implementation, we consider the orthogonal basis on every mesh
edge through the univariate Legendre polynomials and inside every mesh
cell provided by the Gram-Schmidt algorithm applied to the standard
monomial basis.

\medskip
\noindent
Finally, throughout the paper we use the letter $\Cs$ in the error
inequalities to denote a real, positive constant that can have a
different value at any occurrence.
This constant may depend on the model and on some discretization parameters,
such as the coercivity and stability constants of the bilinear form
and of the linear functional used in the variational formulation, the mesh
regularity constants used when defining the properties of the mesh
families to which the numerical method is suitable, and the polynomial
order of the method.
Nevertheless, this constant is always assumed to be independent of the
mesh size parameter $\hh$ that characterizes the mesh and will be
introduced in the next section.\section{The virtual element method}
\label{sec:VEM}

We investigate the performance of the VEM on the elliptic model
problem provided by the Poisson equation with Dirichlet boundary
conditions.
In this section, we briefly review the model equations in strong and
weak form and the formulation of the virtual element approximation.

\medskip
\PGRAPH{The Poisson equation and the virtual element approximation.}
Let $\Omega$ be an open, bounded, connected subset of $\REAL^2$ with polygonal boundary $\Gamma$.
Consider the Poisson equation with homogeneous Dirichlet boundary
conditions in strong from:
\begin{align}
  -\Delta\us = \fs & \qquad\textrm{in~}\Omega,\label{eq:Poisson:strong:A}\\[0.5em]
  \us        = 0   & \qquad\textrm{on~}\Gamma.\label{eq:Poisson:strong:B}
\end{align}
The variational formulation of
problem~\eqref{eq:Poisson:strong:A}-\eqref{eq:Poisson:strong:B} reads
as: \textit{Find $\us\in\HONEzr(\Omega)$ such that}
\begin{align}
  \as(\us,\vs)=\Fs(\vs)\qquad\forall\vs\in\HONEzr(\Omega),
  \label{eq:Poisson:weak}
\end{align}
where the bilinear form
$\as(\cdot,\cdot):\HONE(\Omega)\times\HONE(\Omega)\to\REAL$ is given
by
\begin{align}
  \as(\us,\vs)=\int_{\Omega}\nabla\us\cdot\nabla\vs\dxv
  \label{eq:as:def}
\end{align}
and the right-hand side linear functional $\Fs:\LTWO(\Omega)\to\REAL$
is given by
\begin{align}
  \Fs(\vs)=\int_{\Omega}\fs\vs\dxv,
  \label{eq:Fs:def}
\end{align}
with the (implicit) assumption that $\fs\in\LTWO(\Omega)$.
The well-posedness of the discrete
formulation~\eqref{eq:Poisson:weak} stems from the coercivity and
continuity of the bilinear form $\as(\cdot,\cdot)$, the continuity of
the right-hand side linear functional $\Fs(\cdot)$, and the
application of the Lax-Milgram theorem~\cite[Section 2.7]{Brenner-Scott:2008}.

The numerical method that we consider in this paper is mainly based on
References~\cite{Ahmad-Alsaedi-Brezzi-Marini-Russo:2013,BeiraodaVeiga-Brezzi-Cangiani-Manzini-Marini-Russo:2013},
and provides an optimal approximation on polygonal meshes when the
diffusion coefficient is variable in space.
To ease the presentation, we consider the case of homogeneous
Dirichlet boundary conditions, the extension to the non-homogeneous
case being deemed as straightforward.
Such a case is also considered in the numerical experiments carried
out in this paper.

\medskip
The virtual element approximation of equation~\eqref{eq:Poisson:weak}
reads as: \textit{Find $\ush\in\Vhk$ such that}
\begin{align}
  \ash(\ush,\vsh) = \Fsh(\vsh)
  \qquad\vsh\in\Vhk,
  \label{eq:VEM}
\end{align}
where $\ush$, $\Vhk$, $\ash(\cdot,\cdot)$, $\Fsh(\cdot)$ are the
virtual element approximations of $\us$, $\HONEzr(\Omega)$,
$\as(\cdot,\cdot)$, and $\Fs(\cdot)$.
We review the construction of these mathematical objects in the rest
of this section.

\medskip
\PGRAPH{Mesh notation.}
Let $\calT=\{\Th\}_{\hh\in\calH}$ be a set of decompositions $\Th$ of
the computational domain $\Omega$ into a finite set of nonoverlapping
polygonal elements $\P$.
We refer to $\calT$ as the \emph{mesh family} and to each one of 
its members $\Th$ as the \emph{mesh}.
The subindex label $\hh$, indicating the \textit{mesh size}, is the maximum of the diameters of the mesh
elements, defined by $\hP=\sup_{\xv,\yv\in\P}\ABS{\xv-\yv}$.
We assume that the mesh sizes of the mesh family $\calT$ are in a
countable subset $\calH$ of the real line $(0,+\infty)$ having $0$ as its
unique accumulation point.
Each element $\P$ has a nonintersecting polygonal boundary
$\partial\P$ formed by straight edges $\E$, center of gravity
$\xvP=(\xsP,\ysP)$ and area $\mP$.
We denote the edge mid-point $\xvE=(\xsE,\ysE)$ and its lenght
$\mE$, and with a small abuse of notation, we write $\E\in\partial\P$
to indicate that edge $\E$ is running throughout the set of edges
forming the elemental boundary $\partial\P$.
The convergence analysis of the VEM and the derivation of the error
estimates in the $\LTWO$ and $\HONE$ norms require a few suitable
assumptions on the mesh family $\calT$.
Such assumptions are discussed in detail in the next section.
On every mesh $\Th$, given an integer $k\geq0$, we define the space of
piecewise discontinuous polynomials of degree $k$, $\PS{k}(\Th)$,
containing the functions $\qs$ such that $\restrict{\qs}{\P}\in\PS{k}(\P)$ for
every $\P\in\Th$.

\medskip
\PGRAPH{The virtual element spaces.}
Let $k\geq1$ be an integer number and $\P\in\Th$ a generic mesh
element.
The \emph{conforming virtual element space $\Vhk$ of order $k$ built
on mesh $\Th$} is obtained by gluing together the elemental
approximation spaces denoted by $\Vhk(\P)$, so that
\begin{align}
  \Vhk:=\Big\{\,\vsh\in\HONEzr(\Omega)\,:\,\restrict{\vsh}{\P}\in\Vhk(\P)
  \,\,\,\forall\P\in\Th\,\Big\}.
  \label{eq:Vhk:def}
\end{align}
The local virtual element space $\Vhk(\P)$ is defined in accordance
with the \emph{enhancement strategy} introduced
in~\cite{Ahmad-Alsaedi-Brezzi-Marini-Russo:2013}:
\begin{multline}
  \Vhk(\P) = \bigg\{\,
  \vsh\in\HONE(\P)\,:\,
  \restrict{\vsh}{\partial\P}\in\CS{0}(\partial\P),\,
  \restrict{\vsh}{\E}\in\PS{k}(\E)\,\forall\E\in\partial\P,\,
  \Delta\vsh\in\PS{k}(\P),\,
  \textrm{and~}\\[0.25em]
  \int_{\P}(\vsh-\PinP{k}\vsh)\,\qs\dV=0
  \,\,\forall\qs\in\PS{k}(\P)\backslash\PS{k-2}(\P)
  \,\bigg\},
  \label{eq:VhkP:def}
\end{multline}
where $\PinP{k}$ is the elliptic projection that will be discussed in
the next section; $\PS{k}(\P)$ and $\PS{k}(\E)$ are the linear spaces
of the polynomials of degree at most $k$, which are respectively
defined over an element $\P$ or an edge $\E$ according to our
notation; and $\PS{k}(\P)\backslash\PS{k-2}(\P)$ is the space of
polynomials of degree equal to $k-1$ and $k$.
By definition, the space $\Vhk(\P)$ contains $\PS{k}(\P)$ and the global
space $\Vhk$ is a conforming subspace of $\HONE(\Omega)$.

\medskip
\PGRAPH{The elliptic projection operators.}
The definition in~\eqref{eq:VhkP:def} requires the \emph{elliptic
projection operator} $\PinP{k}:\HONE(\P)\to\PS{k}(\P)$, which, for any
$\vsh\in\Vhk(\P)$, is given by:
\begin{align}  
  \int_{\P}\nabla\PinP{k}\vsh\cdot\nabla\qs\dV  &= \int_{\P}\nabla\vsh\cdot\nabla\qs\dV\quad\forall\qs\in\PS{k}(\P),\label{eq:PiFn}\\[0.5em]
  \int_{\partial\P}\big(\PinP{k}\vsh-\vsh\big)\dS &= 0.                                                              \label{eq:def:Pib_k}
\end{align}
Equation~\eqref{eq:def:Pib_k} allows us to remove the kernel of the
gradient operator from the definition of $\PinP{k}$, so that the
$k$-degree polynomial $\PinP{k}\vsh$ is uniquely defined for every
virtual element function $\vsh\in\Vhk(\P)$.
Moreover, projector $\PinP{k}$ is a polynomial-preserving operator,
i.e., $\PinP{k}\qs=\qs$ for every $\qs\in\PS{k}(\P)$.
We can also define a global projection operator
$\Pin{k}:\HONE(\Omega)\to\PS{k}(\Th)$, which is such that
$\restrict{\Pin{k}\vsh}{\P}=\PinP{k}(\restrict{\vsh}{\P}) \ \forall \P \in \Th$.
A major property of the elliptic projection operator is that every
projection $\PinP{k}\vsh$ of a virtual element function
$\vsh\in\Vhk(\P)$ is computable from the degrees of freedom of $\vsh$
associated with element $\P$ that are defined as
follows.

\medskip
\PGRAPH{The degrees of freedom.}
The degrees of freedom of a virtual element function
$\vsh\in\Vhk(\P)$ are given by the following set of
values~\cite{BeiraodaVeiga-Brezzi-Cangiani-Manzini-Marini-Russo:2013}:
\begin{description}
\item[]\textbf{(D1)} for $k\geq1$, the values of $\vsh$ at the
  vertices of $\P$;

  \medskip
\item[]\textbf{(D2)} for $k\geq2$, the values of $\vsh$ at the $k-1$
  internal points of the $(k+1)$-point Gauss-Lobatto quadrature rule
  on every edge $\E\in\partial\P$.

  \medskip
\item[]\textbf{(D3)} for $k\geq2$, the cell moments of $\vsh$ of order
  up to $k-2$ on element $\P$:
  \begin{align}
    \frac{1}{\mP}\int_{\P}\vsh\,\qs\dV
    \quad\forall\qs\in\PS{k-2}(\P).
    \label{eq:dofs:D3}
  \end{align}
\end{description}
These set of values are unisolvent in $\Vhk(\P)$,
cf.~\cite{BeiraodaVeiga-Brezzi-Cangiani-Manzini-Marini-Russo:2013};
hence, every virtual element function is uniquely identified by them.
The degrees of freedom of a virtual element function in the global
space $\Vhk$ are given by collecting the elemental degrees of freedom
\textbf{(D1)}-\textbf{(D3)}.
Their unisolvence in $\Vhk$ is an immediate consequence of their
unisolvence in every elemental space $\Vhk(\P)$.

\medskip
\PGRAPH{Orthogonal projections.}
From the degrees of freedom of a virtual element function
$\vsh\in\Vhk(\P)$ we can also compute the orthogonal projections
$\PizP{k}\vsh\in\PS{k}(\P)$,
cf.~\cite{Ahmad-Alsaedi-Brezzi-Marini-Russo:2013}.
In fact, the orthogonal projection $\PizP{k}\vsh$ of a function
$\vsh\in\Vhk(\P)$ is the solution of the variational problem:
\begin{align}
  \int_{\P}\Piz{k}\vsh\,\qs\dV =
  \int_{\P}\vsh\,\qs\dV\qquad\forall\qs\in\PS{k}(\P).
\end{align}
The right-hand side is the integral of $\vsh$ against the polynomial
$\qs$, and is computable from the degrees of freedom \textbf{(D3)} of
$\vsh$ when $\qs$ is a polynomial of degree up to $k-2$, and from the
moments of $\PinP{k}\vsh$ when $\qs$ is a polynomial of degree $k-1$
and $k$, cf.~\eqref{eq:VhkP:def}.
Clearly, the orthogonal projection $\PizP{k-1}\vsh$ is also
computable.
As we have done for the elliptic projection, we can also define a
global projection operator $\Piz{k}:\LTWO(\Omega)\to\PS(\Th)$, which
projects the virtual element functions on the space of discontinuous
polynomials of degree at most $k$ built on mesh $\Th$.
This operator is given by taking the elemental $\LTWO$-orthogonal
projection $\PizP{k}\vsh$ in every mesh element $\P$, so that
$\restrict{\big(\Piz{k}\vsh\big)}{\P}=\PizP{k}(\restrict{\vsh}{\P})$,
which is computable from the degrees of freedom of $\vsh$ associated
with element $\P$.

\medskip
\PGRAPH{Approximation properties in the virtual element space.}
Under a suitable regularity assumption on the mesh family used in the 
formulation of the VEM (assumption \textbf{G\ref{g1}} that will be the topic of 
the next section), we can prove the following estimates on the projection and interpolation 
operators:
\begin{enumerate}
\item for every $s$ with $1\le\ss\le\ks+1$ and for every $\ws\in\HS{s}(E)$ there exists a
$\ws_{\pi}\in\PS{k}(\P)$ such that
\begin{equation}
    \ABS{\ws-\ws_{\pi}}_{0,\P} +
    \hP\ABS{\ws-\ws_{\pi}}_{1,E}\leq\Cs\hh^s_{\P}\ABS{\ws}_{s,\P};
\end{equation}
\item for every $s$ with $2\le\ss\leq\ks+1$, for every $\hh$, for all $\P\in\Th$ and for
every $\ws\in\HS{s}(E)$ there exists a $\wsI\in\Vhk(\P)$ such that
\begin{equation}
    \ABS{\ws-\wsI}_{0,\P} +
    \hP\ABS{\ws-\wsI}_{1,\P}\leq\Cs\hP^s\ABS{\ws}_{s,\P}.
\end{equation}  
\end{enumerate}
In these inequalities, $C$ is a real positive constant depending only on the polynomial degree 
$k$ and on some mesh regularity constants that we will introduce and discuss in the next
section.

\medskip
\PGRAPH{The virtual element bilinear forms.}
The elliptic and orthogonal projections are needed to define the
virtual element bilinear form
$\ash(\cdot,\cdot):\Vhk\times\Vhk\to\REAL$, and the forcing term
$\Fsh:\Vhk\to\REAL$.
Following the ``VEM gospel'', we write the discrete bilinear form
$\ash(\cdot,\cdot)$ as the sum of elemental contributions
\begin{align}
  \ash(\ush,\vsh) = \sum_{\P\in\Th}\ashP(\ush,\vsh),\quad
\end{align}
where every elemental contribution is a bilinear form
$\ashP(\cdot,\cdot):\Vhk(\P)\times\Vhk(\P)\to\REAL$ designed to
approximate the corresponding elemental bilinear form
$\asP(\cdot,\cdot):\HONE(\P)\times\HONE(\P)\to\REAL$,
\begin{align*}
  \asP(\vs,\ws) = \int_{\P}\nabla\vs\cdot\nabla\ws\dV,
  \quad\forall\vs,\ws\in\HONE(\P).
\end{align*}
The bilinear form $\ashP(\cdot,\cdot)$ on each element $\P$ is given
by
\begin{align}
  \ashP(\ush,\vsh) 
  = \int_{\P}\nabla\PinP{k}\ush\cdot\nabla\PinP{k}\vsh\dV
  + \SPh\Big( \big(I-\PinP{k}\big)\ush, \big(I-\PinP{k}\big)\vsh \Big).
  \label{eq:ashP:def}
\end{align}
The bilinear form $\SPh(\cdot,\cdot)$ in the definition of
$\ashP(\cdot,\cdot)$ provides the stability term and can be any
computable, symmetric, positive definite bilinear form defined on
$\Vhk(\P)$ for which there exist two positive constants $\cbot$ and
$\ctop$ such that
\begin{align}
  \cbot\asP(\vsh,\vsh)
  \leq\SPh(\vsh,\vsh)
  \leq\ctop\asP(\vsh,\vsh)
  \quad\forall\vsh\in\Vhk(\P)\cap\KER\big(\PinP{k}\big).
  \label{eq:SP:stability}
\end{align}
The inequalities in~\eqref{eq:SP:stability} implies that
$\SPh(\cdot,\cdot)$ scales like $\asP(\cdot,\cdot)$ with respect to
$\hP$.
Also, the stabilization term in the definition of $\ashP(\cdot,\cdot)$
is zero if at least one of its two entries is a polynomial of degree
(at most) $k$, since $\PinP{k}$ is a polynomial preserving operator.

\medskip
In our implementation of the VEM, we consider the stabilization
proposed in~\cite{Mascotto:2018}:
\begin{align}
  \SPh(\vsh,\wsh) =
  \sum_{i=1}^{\NDOFS}\mathcal{A}^{\P}_{ii}\textrm{DOF}_i(\vsh)\textrm{DOF}_i(\wsh),
  \label{eq:stab:D-recipe}
\end{align}
where $\mathcal{A}^{\P}=\big(\mathcal{A}^{\P}_{ij}\big)$ is the matrix
resulting from the implementation of the first term in the bilinear
form $\ashP(\cdot,\cdot)$.
Let $\phi_i$ be the $i$-th ``canonical'' basis functions generating
the virtual element space, which is the function in $\Vhk(\P)$ whose
$i$-th degree of freedom for $i=1,\ldots,\NDOFS$ (according to a
suitable renumbering of the degrees of freedom in \textbf{(D1)},
\textbf{(D2)}, and \textbf{(D3)}), has value equal to $1$ and all
other degrees of freedom are zero.
These basis function are unknown in the virtual element framework, but
their projections $\PizP{k-1}\nabla\phi_i$ (and
$\PizP{k-1}\nabla\phi_j$) are computable from their degrees of
freedom.
With this notation, the $i,j$-th entry of matrix $\mathcal{A}^{\P}$ is
given by
\begin{align}
  \mathcal{A}^{\P}_{ij}:=\asP\big(\PinP{k}\phi_i,\PinP{k}\phi_j\big).
\end{align}
The stabilization in~\eqref{eq:stab:D-recipe} is sometimes called the
``\emph{D-recipe stabilization}'' in the virtual element literature,
and contains the so called ``\emph{dofi-dofi (dd)
stabilization}'' originally proposed
in~\cite{BeiraodaVeiga-Brezzi-Cangiani-Manzini-Marini-Russo:2013} as the special case with $\mathcal{A}_{ii}=1$:
\begin{align}
  \SPhdd(\vsh,\wsh) =
  \sum_{i=1}^{\NDOFS}\textrm{DOF}_i(\vsh)\textrm{DOF}_i(\wsh).
  \label{eq:stab:dofi-dofi}
\end{align}
We explicitly mention the stabilization~\eqref{eq:stab:dofi-dofi}
because many convergence results available from the literature, which
we briefly review in Appendix~\ref{app:error}, are obtained by using it.

\medskip
The stabilization term, and,
in particular, condition~\eqref{eq:SP:stability}, is designed so
that $\ashP(\cdot,\cdot)$ satisfies the two fundamental properties:
\begin{description}
\item[-] {\emph{$k$-consistency}}: for all $\vsh\in\Vhk(\P)$ and for all
  $\qs\in\PS{k}(\P)$ it holds that
  \begin{align}
    \label{eq:k-consistency}
    \ashP(\vsh,\qs) = \asP(\vsh,\qs);
  \end{align}
  
  \medskip
\item[-] {\emph{stability}}: there exist two positive constants
  $\alpha_*,\,\alpha^*$, independent of $\hh$ and $\P$, such that
  \begin{align}
    \label{eq:stability}
    \alpha_*\asP(\vsh,\vsh)
    \leq\ashP(\vsh,\vsh)
    \leq\alpha^*\asP(\vsh,\vsh)\quad\forall\vsh\in\Vhk(\P).
  \end{align}
\end{description} 

\medskip
\PGRAPH{The virtual element forcing term.}
To approximate the right-hand side of~\eqref{eq:VEM}, we split it into
the sum of elemental contributions and every local linear functional is
approximated by using the orthogonal projection $\PizP{k}\vsh$:
\begin{align}
  \Fs(\vsh) = \sum_{\P\in\Th}\big(\fs,\PizP{k}\vsh\big)_{\P}.
  \textrm{~~where~~}\big(\fs,\PizP{k}\vsh\big)_{\P} = \int_{\P}\fs\,\PizP{k}\vsh\dV.
  \label{eq:extensionc}
\end{align}

\medskip
\PGRAPH{Main convergence properties.}
The well-posedness of the discrete formulation~\eqref{eq:VEM} stems
from the coercivity of the bilinear form $\ash(\cdot,\cdot)$, the
continuity of the right-hand side linear functional
$\big(\fs,\Piz{k}\cdot\big)$ and the application of the Lax-Milgram
theorem~\cite[Section 2.7]{Brenner-Scott:2008}.

\medskip
In this work, we are interested in checking whether the VEM mantains
optimal convergence rates on different mesh families that may display
some pathological situations.
From a theoretical viewpoint, the convergence estimates hold under
some constraints on the shapes of the elements forming the
mesh, called \emph{mesh geometrical (or regularity) assumptions}.
We summarize the major findings from the literature in
Appendix~\ref{app:error} and in the next sections we will
investigate how breaking such constraints may affect these results.

\medskip
Let $\us\in\HS{k+1}(\Omega)$ be the solution to the variational
problem~\eqref{eq:Poisson:weak} on a convex domain $\Omega$ with
$\fs\in\HS{k}(\Omega)$.
Let $\ush\in\Vhk$ be the solution of the virtual element method
\eqref{eq:VEM} on every mesh of a mesh family $\calT=\{\Th\}$
satisfying a suitable set of mesh geometrical assumptions.
Then, a strictly positive constant $\Cs$ exists such that
\begin{itemize}
\item the $\HONE$-error estimate holds:
  \begin{align}
    \label{eq:source:problem:H1:error:bound}
    \norm{\us-\ush }{1,\Omega}\leq
    \Cs\hh^{k}\left(
    \norm{\us}{k+1,\Omega} 
    + \snorm{\fs}{k,\Omega}
    \right);
  \end{align}

  \medskip
\item the $\LTWO$-error estimate holds:
  \begin{align}
    \label{eq:source:problem:L2:error:bound}
    \norm{\us-\ush}{0,\Omega}\leq 
    \Cs\hh^{k+1}\left(
    \norm{\us}{k+1,\Omega}
    + \snorm{\fs}{k,\Omega}
    \right).
  \end{align}
\end{itemize}
Constant $\Cs$ may depend on the stability constants $\alpha_*$ and
$\alpha^*$, on mesh regularity constants which we will
introduce in the next section, on the size of the computational domain
$\ABS{\Omega}$, and on the approximation degree $k$.
Constant $\Cs$ is normally independent of $\hh$, but for the most
extreme meshes it may depend on the ratio between the longest and
shortest edge lenghts, cf.  Appendix~\ref{app:error}.

\medskip
Finally, we note that the approximate solution $\ush$ is not
explicitly known inside the elements.
Consequently, in the numerical experiments of
Section~\ref{subsec:performance}, we approximate the error
in the $\LTWO$-norm as follows:
\begin{align*}
  \norm{\us-\ush}{0,\Omega}\approx
  \norm{\us-\Piz{k}\ush}{0,\Omega},
\end{align*}
where $\Piz{k}\ush$ is the global $\LTWO$-orthogonal projection of the
virtual element approximation $\ush$ to $\us$.
On its turn, we approximate the error in the energy norm as follows:
\begin{align*}
  \snorm{ \us-\ush }{1,\Omega}\approx\ash(\usI-\ush,\usI-\ush)
\end{align*}
where $\usI$ is the virtual element interpolant of the exact solution
$\us$.\section{Geometrical Assumptions}
\label{sec:geometrical:assumptions}
In this section, we review the geometrical assumptions appeared in the
VEM literature since their definition in
\cite{BeiraodaVeiga-Brezzi-Cangiani-Manzini-Marini-Russo:2013}.
All the assumptions are defined for a single mesh $\Th$, 
but the conditions contained in them are required to hold 
independently of $h$.
Therefore, when considering a mesh family $\calT=\{\Th\}_h$, 
these assumptions have to be verified simultaneously by every
$\Th \in \calT$.\\
It is well-known from the FEM literature that the approximation
properties depend on specific assumptions on the geometry 
of the elements.
For example, classical geometrical assumptions for a family of
triangulations $(\Th)_{\hh\to 0}$, 
are the ones respectively introduced in \cite{ciarlet2002finite} 
and \cite{zlamal1968finite}:

\smallskip
\begin{description}
    \item [$(a)$] \textit{Shape regularity condition:} there exists a real number $\gamma\in(0,1)$, independent of $h$, such that $\forall \P \in \Omega_h$ we have
    \[\gamma h_\P\le r_\P,\]
    where $h_\P$ and $r_\P$ are, respectively, the longest edge in $\P$ and its inradius;
    \item [$(b)$] \textit{Minimum angle condition:} there exists $\alpha_0>0$, independent of $h$, such that $\forall \P \in \Omega_h$ we have
    \[\alpha_\P\ge\alpha_0,\]
    where $\alpha_\P$ is the minimal angle of $\P$.
\end{description}

\medskip
Similarly, in the VEM we need a set of geometrical assumptions to ensure approximation properties.
The first pair of assumptions were proposed in
\cite{BeiraodaVeiga-Brezzi-Cangiani-Manzini-Marini-Russo:2013} 
and remained untouched 
also in~\cite{Ahmad-Alsaedi-Brezzi-Marini-Russo:2013} 
and \cite{Brenner:2017:SEV}.
In these papers, the Authors assume that a real constant $\rho\in(0,1)$
exists, independent of $\hh$, such that two conditions hold:

\smallskip
\begin{gass}
  Every polygonal cell $\P\in\Th$ is star-shaped with respect 
  to a disc with radius $\rho\hP$.
  \label{g1}
\end{gass}

\begin{gass}
  For every polygonal cell $\P\in\Th$, the length $\mE$ of every edge $\E\in\partial\P$ satisfies
  $\mE\ge\rho\hP$.
  \label{g2}
\end{gass}

\smallskip
\noindent
Constant $\rho$ is often referred to as \textit{mesh regularity 
constant} or \textit{parameter}.
Condition \textbf{G\ref{g1}} can be weakened in
the following way, as specified in
\cite{BeiraodaVeiga-Brezzi-Cangiani-Manzini-Marini-Russo:2013} and
more accurately in \cite{Brenner:2017:SEV}:

\smallskip
\noindent
\textbf{Assumption G1 - weak} \ \textit{Every polygonal cell
  $\P\in\Th$ is the union of a finite number $\Ns$ of disjoint
  polygonal subcells
  $\P_1,\dots,\P_{\Ns}$ such that, for $j=1,\dots,\Ns$,
  \begin{description}
  \item [$(a)$] element $\P_j$ is star-shaped with respect to a disc 
  with radius $\rho{\hP}_j$;
  \item [$(b)$] elements $\P_j$ and $\P_{j+1}$  share a common edge.
\end{description}}

\smallskip
Assumption \textbf{G\ref{g1}} (or \textbf{G1 - weak}) is the
polygonal extension of the classical conditions for triangular
meshes, with $\hP$ indicating the elemental diameter instead 
of the longest edge.
Under assumption \textbf{G1 - weak}, and therefore also under
\textbf{G\ref{g1}}, it can be proved \cite{Brenner:2017:SEV} 
that the simplicial triangulation of $\P$ determined by the
star-centers (the centers of the discs in \textbf{G\ref{g1}} and
\textbf{G1 - weak}) of $\P_1,\dots,\P_{\Ns}$ satisfies the
\textit{shape regularity} and the \textit{minimum angle} conditions.
Moreover, for $1\le j,k\le N$ it holds that
${h_{\P_{j}}}\slash{h_{\P_{k}}}\le\rho^{-|j-k|}$.\\
These assumptions are more restrictive than necessary, but at 
the same time they are not particularly demanding, since they 
allow the method to work on very general decompositions.
This fact was already mentioned in the very first papers.
For example, in~\cite[Ahmad et al.]{Ahmad-Alsaedi-Brezzi-Marini-Russo:2013} 
the Authors say that: 
\begin{quote}
  \textit{``Actually, we could get away with even more general
    assumptions, but then it would be long and boring to make precise
    (among many possible crazy decompositions that nobody will ever
    use) the ones that are allowed and the ones that are not.''}
\end{quote} 

\medskip
In \cite{Beirao-Lovadina-Russo:2016} and \cite{brenner2018virtual}
assumption~\textbf{G\ref{g1}} is preserved, but
assumption~\textbf{G\ref{g2}} is substituted by the alternative version:

\smallskip
\begin{gass}
  There exists a positive integer $\Ns$, independent of $\hh$, such
  that the number of edges of every polygonal cell $\P\in\Th$ is (uniformly)
  bounded by $\Ns$.
  \label{g3}
\end{gass}

\smallskip
\noindent
Assumption~\textbf{G\ref{g2}} implies assumption~\textbf{G\ref{g3}}.
However, assumption~\textbf{G\ref{g3}} is weaker than
assumption~\textbf{G\ref{g2}}, as it allows for edges arbitrarily
small with respect to the element diameter.
Both assumption pairs~\textbf{G\ref{g1}}+\textbf{G\ref{g2}} 
%% and assumption
and~\textbf{G\ref{g1}}+\textbf{G\ref{g3}} imply that the number
of vertices of $\P$ and the minimum angle of the simplicial
triangulation of $\P$ given by connecting the vertices of $\P$ and its
star-center, are controlled by $\rho$.

\begin{figure}[htbp]
  \centering
	\includegraphics[width=.6\linewidth]{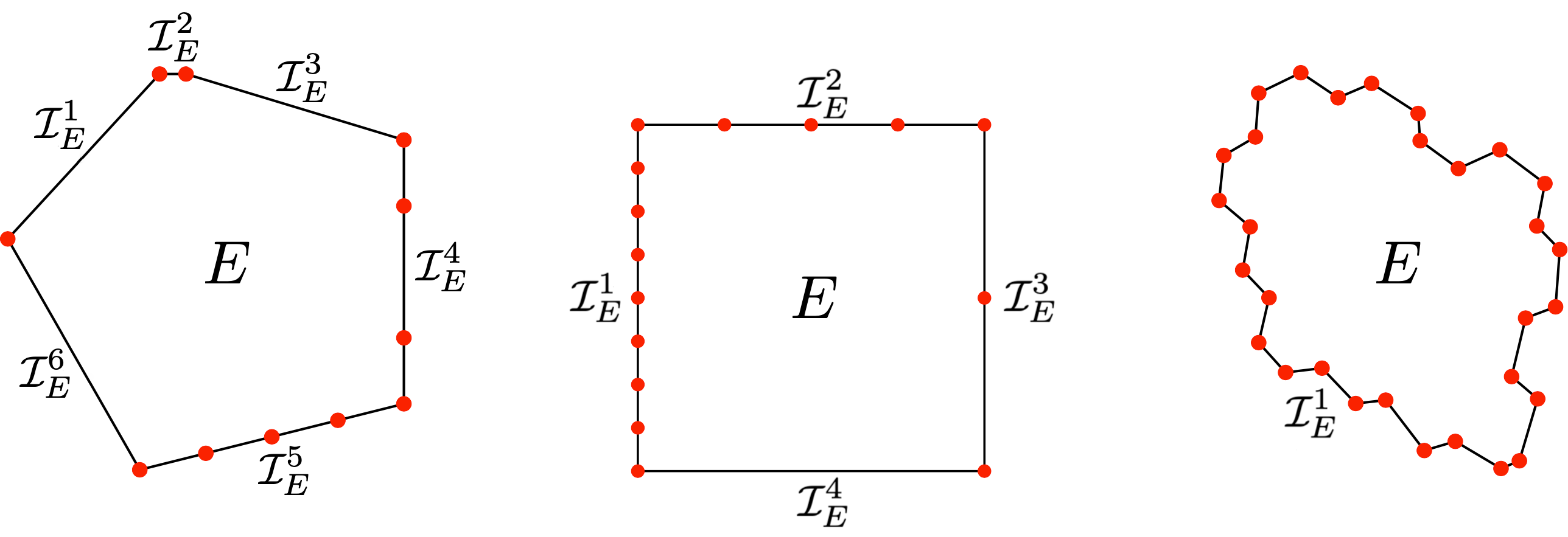} %% {Figures/assumption-g4.png}
  \caption{Examples of admissible elements according to assumption~\textbf{G\ref{g4}}. Red dots indicate the vertices of the element.}
  \label{fig:g4}
\end{figure}

\medskip
Another step forward in the direction of refining the geometrical
assumptions has been made in \cite{da2020sharper}.
In addition to assumption~\textbf{G\ref{g1}}, the Authors imagine to
unwrap the boundary $\partial\P$ of each polygon $\P \in \Omega_h$
onto an interval $I_\P$ of the real line, obtaining a
one-dimensional mesh $\calI_\P$.
The collection of the unwrapped boundaries of all elements in a mesh $\Th$ is
denoted by $\{\calI_\P\}_{\P\in\Th}$.
Moreover, each one-dimensional mesh $\calI_\P$ can be 
subdivided into a number of disjoint sub-meshes $\calI_\P^1, \ldots, \calI_\P^N$, corresponding to the edges of $\P$ (we consider each $\calI_\P^i$ as a mesh as it may contain more than one edge, see Fig.~\ref{fig:g4}).
Then, the following condition is assumed.

\smallskip
\begin{gass} 
  \label{g4}
  For every polygonal cell $\P\in\Th$, the family $\{\calI_\P\}_{\P \in
    \Omega_h}$ is piecewise quasi-uniform, that is:
  \begin{description}
  \item [$(a)$] each mesh $\calI_\P$ can be subdivided into at most
    $N$ disjoint sub-meshes
    $\calI_\P^1, \ldots, \calI_\P^N$, for some
    $N\in \mathbb N$;
  \item [$(b)$] each sub-mesh $\calI_\P^i$, $i=1, \ldots, N$, is quasi
    uniform: the ratio between the largest and the smallest element in
    $\calI_\P^i$ is bounded from above by some $\cs\in\mathbb R^+$ independent of $h$.
  \end{description}
\end{gass}

\smallskip
Each polygon $\P$ is in a one-to-one correspondence to a one-dimensional
mesh $\calI_\P$, but a sub-mesh $\calI_\P^i\subset\calI_\P$ might
contain more than one edge of $\P$.
This implies that assumption~\textbf{G\ref{g4}} does not require a
uniform bound on the number of edges in each element and does not
exclude the presence of small edges, cf. Fig.~\ref{fig:g4}.
For instance, the mesh families created by agglomeration, cracking, gluing, etc.. of
existing meshes are admissible according to~\textbf{G\ref{g4}}.
\\
According to the literature (see Appendix~\ref{app:error}), possible 
assumption pairs requested to guarantee the convergence of the VEM are given by 
combining \textbf{G\ref{g1}} (or, equivalently, \textbf{G1 - weak}) with either 
\textbf{G\ref{g2}} or
\textbf{G\ref{g3}} or \textbf{G\ref{g4}}.
\section{Breaking the geometrical assumptions}
\label{sec:breaking}
In this section, we test the behaviour of the virtual element method
on a number of mesh ``datasets'', to stress one or more of the 
geometrical assumptions discussed in Section~\ref{sec:geometrical:assumptions}.
We call a \emph{dataset} a collection $\calD:=\{\Omega_n\}_{n=0,\ldots,N}$ 
of meshes $\Omega_n$ covering the domain $\Omega=(0,1)^2$ such that 
\begin{description}
    \item [-] the mesh $\Omega_{n+1}$ has smaller mesh size than
    $\Omega_{n}$ for every $n=0,\ldots,N-1$;
    \item [-] the meshes $\Omega_n$ follow a common refinement pattern, 
    so that they contain similar polygons organized in similar configurations.
\end{description}
Note that each mesh $\Omega_n$ is uniquely identified by its size
as $\Th$, therefore we can consider a dataset $\calD$ as a subset 
of a mesh family: $\calD=\{\Th\}_{h\in\calH'}\subset\calT$ 
where $\calH'$ is a finite subset of $\calH$.

\medskip
In addition to the violation of the geometrical assumptions, we are 
also interested in the behaviour of the VEM when the
measures of mesh elements and edges scale in a nonuniform way in the 
refinement process.
To this end,
for each mesh $\Omega_n\in\calD$ 
we define the following quantities and study their trend for $n\to N$:
\begin{align}
A_n=\frac{\max_{\P\in \Omega_n} |\P|}{\min_{\P\in \Omega_n} |\P|}
\quad\textrm{and}\quad 
e_n=\frac{\max_{e\in \Omega_n} |e|}{\min_{e\in \Omega_n} |e|}.
\label{eq:scaling:indicators}
\end{align}
We specifically designed six datasets in order to consider several
possible combinations of the geometrical assumptions 
of the previous section and the scaling indicators $A_n$ and $e_n$,
as shown in~Table~\ref{table:ratios}.
Note that most of the considered datasets do not fulfill any set of 
geometrical assumptions required by the convergence analysis found 
in the literature (see Appendix~\ref{app:error}).

\subsection{Datasets definition}
\label{subsec:datasets}
We now introduce the datasets, describing for each of them how they are built, which geometrical assumptions they
fulfill or violate, and how the indicators $A_n$ and $e_n$ depend on $n$ in the limit for $n\to\Ns$.
Each dataset is built around (and often named after) a particular polygonal element contained in it, which is meant to stress one or more assumptions or indicators.
The detailed construction algorithms, together with the explicit computations of $A_n$ and $e_n$ for all datasets, can be found in Appendix~\ref{app:dataset}.

\medskip
\PGRAPH{Reference dataset.}
The first dataset, $\Dtriangle$, contains only triangular meshes that 
are built by inserting a number of vertices in the domain through the Poisson Disk Sampling algorithm \cite{bridson2007fast}, and connecting them in a Delaunay triangulation (see Appendix~\ref{app:dataset:reference}). 
The refinement is obtained by increasing the number of vertices generated by the Poisson algorithm.
The meshes in this dataset do not violate any of the geometrical assumptions and the indicators $A_n$ and $e_n$ are almost constant.
We use $\Dtriangle$ as the reference dataset to evaluate the other datasets
by comparing the performance of the VEM over them. 

\medskip
\PGRAPH{Hybrid datasets.}
Next, we consider some hybrid datasets, characterized by a progressive insertion in $\Omega$ of one or more identical polygonal elements (called the \emph{initial polygons}), the rest of the domain being tessellated by triangles. 
These triangles are created by the library \textit{Triangle} \cite{Shewchuk:2005},
bounding the area of the triangular elements with the area
of the initial polygons.
Steiner points~\cite{Shewchuk:2005} can be added, and the edges 
of the initial polygons are split when necessary by the insertion of new 
vertices.
The refinement is iterative, with parameters to indicate size, shape and number of the initial polygons; details on this process are provided in the Appendix~\ref{app:dataset:hybrid}.

\begin{figure}[htbp]
\centering
	\includegraphics[width=.7\linewidth]{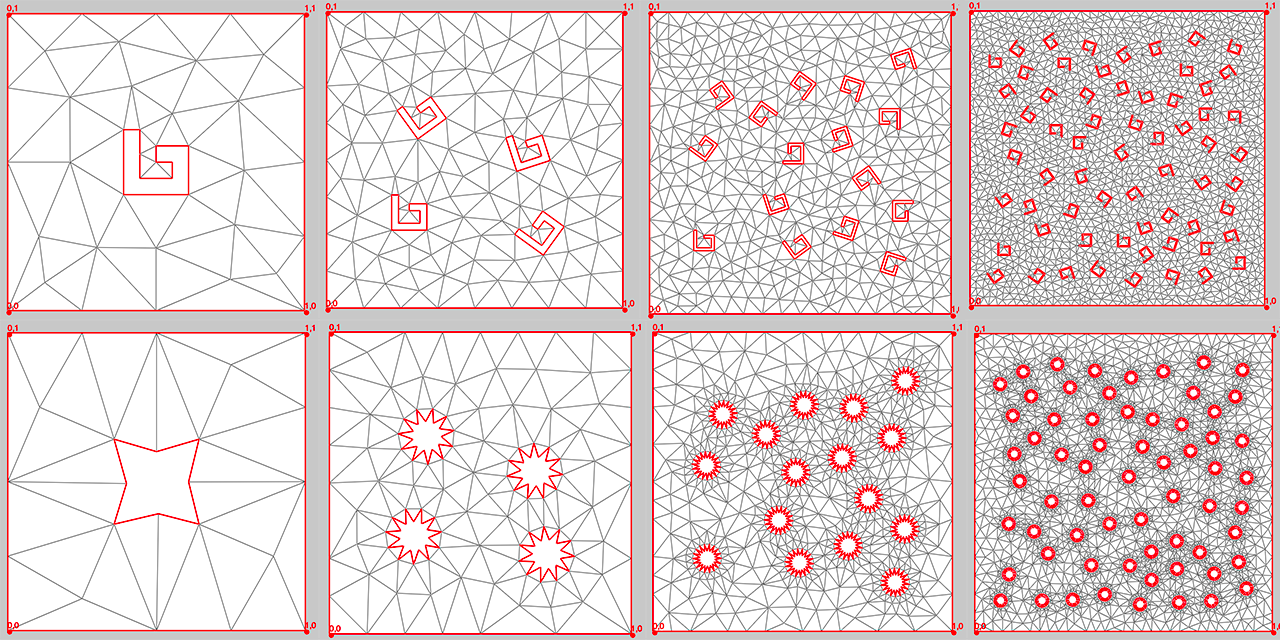} %% {Figures/dataset-hybrid.png}
\caption{Meshes $\Omega_0, \Omega_2, \Omega_4, \Omega_6$ from datasets $\Dmaze$ (top) and $\Dstar$ (bottom).}
\label{fig:hybrid-mesh}
\end{figure}

\medskip
The top and bottom panels of Fig.~\ref{fig:hybrid-mesh} respectively show the 
datasets $\Dmaze$ and $\Dstar$, which we selected as they violate different geometrical assumptions. 
Other choices for the initial polygons are possible, for instance considering the ones
in Benchmark \cite{attene2019benchmark}.

\medskip
A ``maze'' is a $10$-sided polygonal element $\P$ spiralling around an external 
point.
Progressively, each mesh in $\Dmaze$ contains an increasing 
number of mazes $\P$ with decreasing thickness as $n\to N$.
Every $\P$ is obviously \textit{not} star-shaped, challenging assumption~\textbf{G\ref{g1}}.
Moreover, the length of the shortest edge $\E$
of $\P$ decreases faster than the diameter $\hP$ of the polygon.
This fact implies, on one side, that the ratio $|e|/\hP$ of
assumption~\textbf{G\ref{g2}} cannot be bounded from below by a constant $\rho$ that
is independent of $h$, and, on the other side, that assumption~\textbf{G\ref{g1}-weak} also fails.
Indeed, even splitting $\P$ into a finite number of rectangles, it is not possible
to define a global radius $\rho$, independent of $\hh$, with respect to which the union of these
rectangles is star-shaped according to \textbf{G1}, if the shortest edge of $\P$ 
is constantly decreasing.
Concerning the scaling indicators, we have $A_n\sim a^n$ for a constant $e<a<3$ and $e_n\sim n\log(n)$.

\begin{figure}[htbp]
\centering
	\includegraphics[width=.45\linewidth]{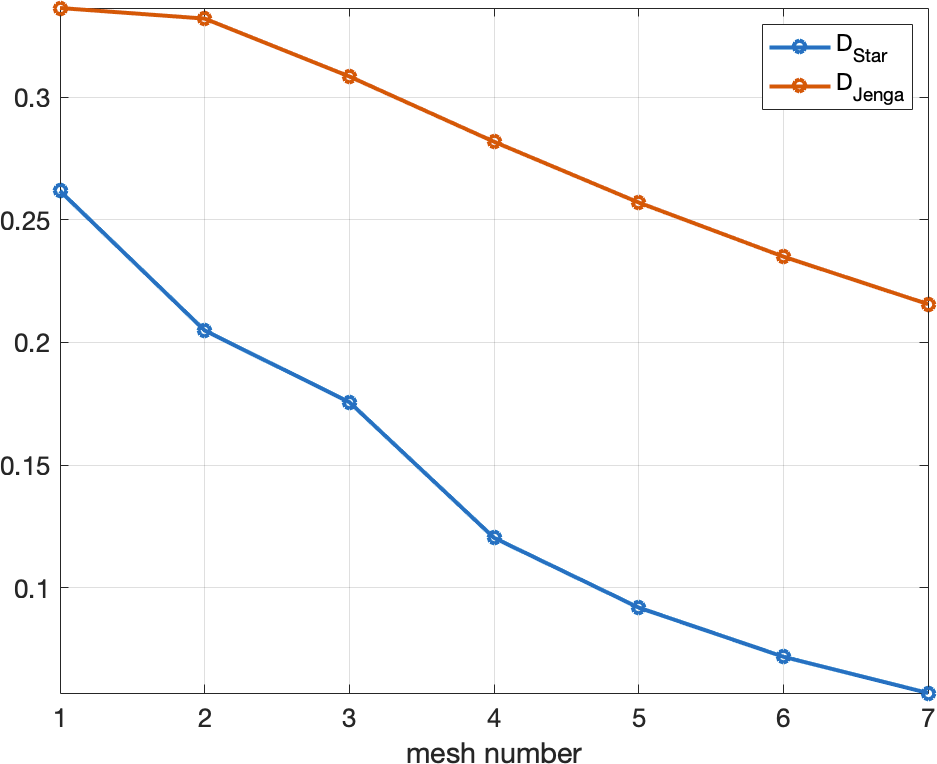} %% {Figures/diameter-radius-ratio.png}
\caption{Ratio $r/\hP$ for datasets $\Dstar$ and $\Djenga$.}
\label{fig:ratio}
\end{figure}

\medskip
Dataset $\Dstar$ is built by inserting star-like polygonal elements,
still denoted by $\P$. 
As $n\to N$, the number of spikes of each $\P$ increases and the inner vertices
of the star move towards the barycenter of the element. 
In this case, assumption~\textbf{G\ref{g3}} is not satisfied because the number 
of spikes in each $\P$ increases from mesh to mesh.
Therefore, the total number of vertices and edges in a single element cannot be bounded uniformly.\\
Last, each star $\P$ is star-shaped with respect to the maximum circle inscribed in it.
However, as shown in Fig.~\ref{fig:ratio}, the radius $r$ of such circle decreases faster than the elemental diameter $\hP$, therefore it is not possible to define a global $\rho>0$ able to uniformly bound from below the quantity $r/\hP$: this violates assumption~\textbf{G\ref{g1}}.
In order to satisfy assumption~\textbf{G1-weak}, we should split each $\P$ into a number of 
sub-polygons that are star-shaped according to \textbf{G\ref{g1}}.
Independently of the way we partition $\P$, the number of sub-polygons would always be
bigger than or equal to the number of spikes in $\P$, which is constantly increasing.
So, the number of sub-polygons would tend to infinity violating
condition \textbf{G1-weak}.
Last, both $A_n$ and $e_n$ scale linearly.

\medskip
\PGRAPH{Mirroring datasets.}
Another possible strategy to build a sequence of meshes whose elements are
progressively smaller, is to adopt a \textit{mirroring} technique.
In practice, we start with the first base mesh $\widehat{\Omega}_0$, which coincides with the first computational mesh $\Omega_0$.
At every step $n\geq1$, we build a new base mesh $\widehat{\Omega}_n$
from the previous base mesh $\widehat{\Omega}_{n-1}$.
The computational mesh $\Omega_{n}$ is then obtained by mirroring $\widehat{\Omega}_{n}$ $4^n$ times
and resizing everything to fit the domain $\Omega$.
This construction allows us to obtain a number of vertices and degrees of freedom
in each mesh that is comparable to that of the meshes at the same refinement level in datasets 
$\Dmaze$ and $\Dstar$.\\
Examples of meshes from mirrored datasets are presented in Fig.~\ref{fig:nonuniform-mesh}; examples of non-mirrored base meshes are visible in Appendix~\ref{app:dataset:mirroring}.
Algorithms for the construction of the following datasets are reported in Appendix~\ref{app:dataset:mirroring}, while the mirroring algorithm is detailed in Appendix~\ref{app:dataset:mirroring-algorithm}.

\medskip
In the case of the dataset $\Djenga$, we build the $n$-th base mesh $\widehat{\Omega}_{n}$ as follows.
We start by drawing two horizontal edges that split the domain $(0,1)^2$ 
into three horizontal rectangles with area equal to $1/4$, $1/2$ and $1/4$ respectively.
Then, we split the rectangle with area $1/2$ vertically, into two equally-sized rectangles with area $1/4$.
This provides us with base mesh $\widehat{\Omega}_0$, which coincides with mesh $\Omega_0$.
At each next refinement step $n\geq1$, we split the left-most rectangle in the middle 
of the base mesh $\widehat{\Omega}_{n-1}$ by adding a new vertical edge, and apply the mirroring 
technique to obtain $\Omega_{n}$.
This process is shown in the top panels of Fig.~\ref{fig:nonuniform-mesh}.

This mesh family breaks all assumptions~\textbf{G\ref{g1}} (and~\textbf{G1-weak}),
\textbf{G\ref{g2}}, \textbf{G\ref{g3}}, and~\textbf{G\ref{g4}}.
In fact, the length of the radius $r$ of the biggest possible disc inscribed into a rectangle is equal to $1/2$ of its shortest edge $\E$.
As shown in Fig.~\ref{fig:ratio}, the ratio $|\E|/h_\P$, decreases unboundedly in the 
left rectangle $\P$ every time we split it, and consequently $r/\hP$ decreases at a similar rate.
This implies that a lower bound with a uniform constant $\rho$ independent of $\hh$ cannot exist for these ratios, thus breaking 
assumptions \textbf{G\ref{g1}}, \textbf{G1-weak} and \textbf{G\ref{g2}}.
In addition, the number of edges of the top and bottom rectangular elements also grows unboundedly, against assumption \textbf{G\ref{g3}}.
Last, the one-dimensional mesh of assumption \textbf{G\ref{g4}}, which is
built on the elemental boundary of the top and bottom rectangular elements,
cannot be subdivided into a finite number of quasi uniform sub-meshes.
In fact, either we have infinite sub-meshes or an infinite edge ratio.
Finally, we note that both $A_n$ and $e_n$ scale like $2^{n}$.

\medskip
In the case of the dataset $\Dslices$ (Fig.~\ref{fig:nonuniform-mesh}, middle), we build the $n$-th 
base mesh $\widehat{\Omega}_{n}$ as follows.
First, we sample a collection of points along the diagonal 
(the one connecting the vertices with coordinates $(0,1)$ and $(1,0)$)
of the reference square $[0,1]^2$, and connect them to the vertices $(0,0)$ and $(1,1)$.
In particular, at each step $n\geq0$, the base mesh $\widehat{\Omega}_{n}$ contains the vertices $(0,0)$ and $(1,1)$, plus the vertices
with coordinates $(2^{-i},1-2^{-i})$ and $(1-2^{-i},2^{-i})$ for $i=1,\ldots,n+2$.
Then, we apply the mirroring technique.

The dataset $\Dslices$ violates assumptions \textbf{G\ref{g1}} and 
\textbf{G1-weak}.
In fact, up to a multiplicative scaling factor depending on $\hh$,
the length of the radius of the biggest inscribed disc in every element $\P$ is decreasing 
faster than the diameter of the element, which is constantly equal to $\sqrt{2}$ times the same 
scaling factor, thus violating \textbf{G\ref{g1}}.
Furthermore, the dataset also breaks assumption \textbf{G1-weak} because any finite subdivisions of its elements
would suffer the same issue.
Instead, the other geometrical assumptions are satisfied.
Since no edge is split, we find that $e_n \sim c$, while $A_n \sim 2^{n}$.

\medskip
In $\Dulike$ (Fig.~\ref{fig:nonuniform-mesh}, bottom), we build $\widehat{\Omega}_{n}$ at each step $n\geq0$ by inserting $2^n$ equispaced $U$-shaped continuous polylines inside the domain, creating as many $U$-like polygons.
Then, we apply the mirroring technique.

For arguments similar to the ones brought for $\Dmaze$, $\Dulike$
does not satisfy assumptions \textbf{G\ref{g1}}, \textbf{G\ref{g1}-weak} 
and \textbf{G\ref{g2}}.
For connectivity reasons, the lower side of the outer $U$-shaped polygon of every base mesh must be split into smaller segments when we apply the mirroring technique.
Therefore, the number of edges of such cells cannot be limited from above, contradicting 
assumption \textbf{G\ref{g3}}.
Nonetheless, assumption \textbf{G\ref{g4}} is satisfied because this subdivision is uniform.
Last, edge lengths scale exponentially and areas scale uniformly, i.e., $e_n \sim 2^{n}$, $A_n\sim c$.
\begin{figure}[htbp]
\centering
	\includegraphics[width=.7\linewidth]{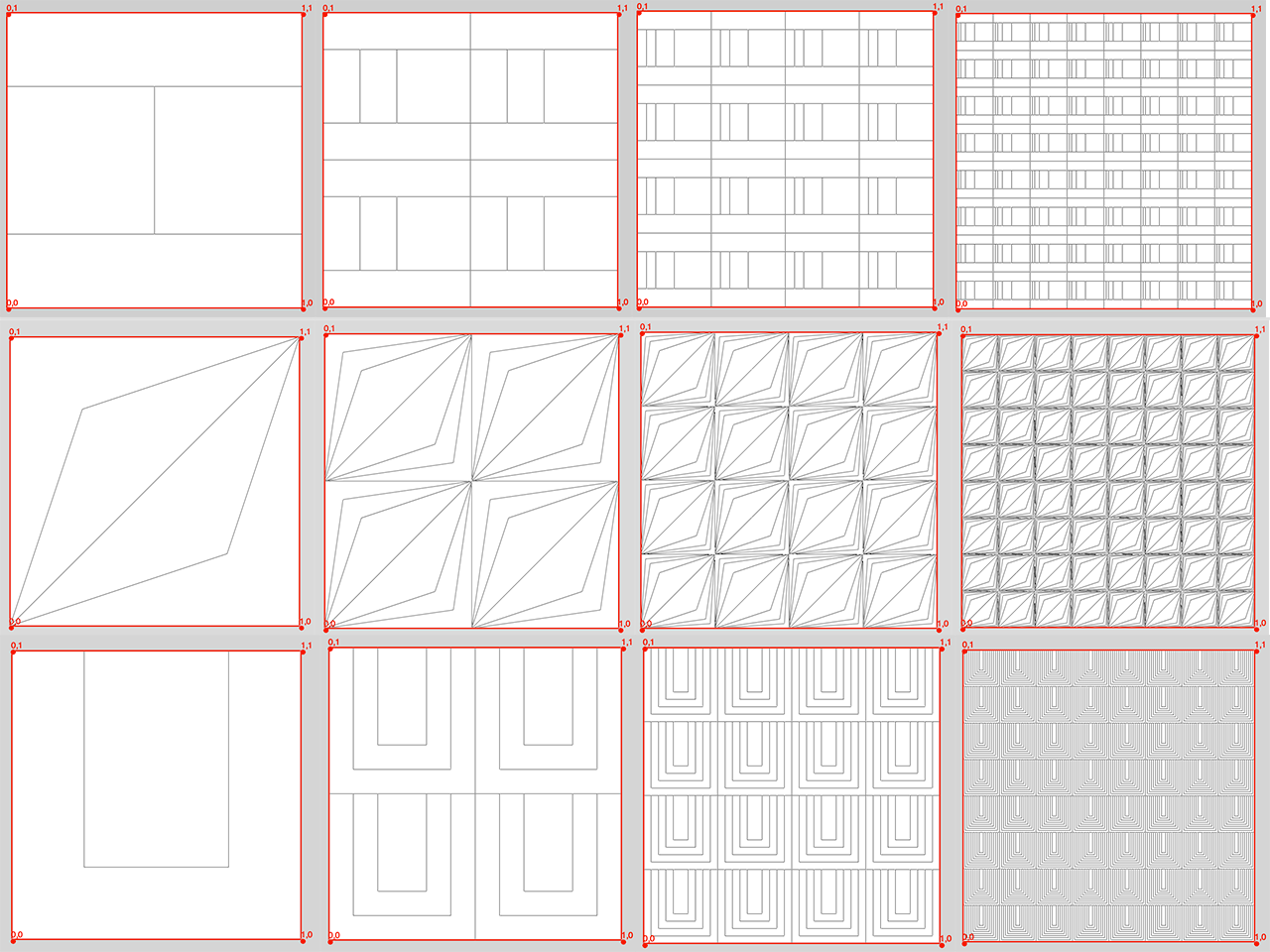} %% {Figures/dataset-mirroring.png}
\caption{Meshes $\Omega_0, \Omega_1, \Omega_2, \Omega_3$ from datasets $\Djenga$ (top), $\Dslices$ (middle) and $\Dulike$ (bottom).}
\label{fig:nonuniform-mesh}
\end{figure}

\medskip
\PGRAPH{Multiple mirroring datasets.}
As a final test, we modified datasets $\Djenga$, $\Dslices$ and $\Dulike$ in order to stress the indicators $A_n$ and $e_n$ harder.\\
This is easily obtained by inserting four new elements at each step instead of one, as explained in Appendix~\ref{app:dataset:multiple-mirroring}.
The resulting datasets, $\DjengaM$, $\DslicesM$ and $\DulikeM$, are qualitatively similar to the mirroring datasets above.
These datasets fulfill the same assumptions as their respective original versions, but the number of elements at each refinement step now increases four times faster.
The indicators $A_n$ and $e_n$ change from $2^n$ to $2^{4n}$, but $A_n$ remains constant for $\DulikeM$, and $e_n$ remains constant for $\DslicesM$.

\begin{table}[t]
\centering
\caption{Summary of the geometrical conditions violated and the asymptotic trend of the indices $A_n$ and $e_n$ for all datasets ($a$ is a constant such that $e<a<3$). Assumption \textbf{G1-weak} is not explicitly reported because all the considered datasets that violate \textbf{G1}, also violate \textbf{G1-weak}.}
\label{table:ratios}
\vspace{0.5cm}
\begin{tabular}{ccccccc}
\hline\noalign{\smallskip}
\textbf{dataset \ } & $\Dtriangle$ & $\Dmaze$ & $\Dstar$ & $\Djenga$ & $\Dslices$ & $\Dulike$\\
\noalign{\smallskip}\hline\noalign{\smallskip}

\textbf{G1} & 
& $\times$ & $\times$ & $\times$ & $\times$ & $\times$ \\

\textbf{G2} & 
& $\times$ & & $\times$ & & $\times$ \\

\textbf{G3} & 
& & $\times$ & $\times$ & & $\times$ \\

\textbf{G4} & 
& & & $\times$ & & \\

$A_n$ & 
$c$ & $a^n$ & $n$ & $2^n$ & $2^n$ & $c$ \\

$e_n$ & 
$c$ & $n\log(n)$ & $n$ & $2^n$ & $c$ & $2^n$\\

\noalign{\smallskip}\hline
\end{tabular}
\end{table}

\subsection{Performance analysis}
\label{subsec:performance}
\begin{figure}[htbp]
\centering
	\includegraphics[width=\linewidth]{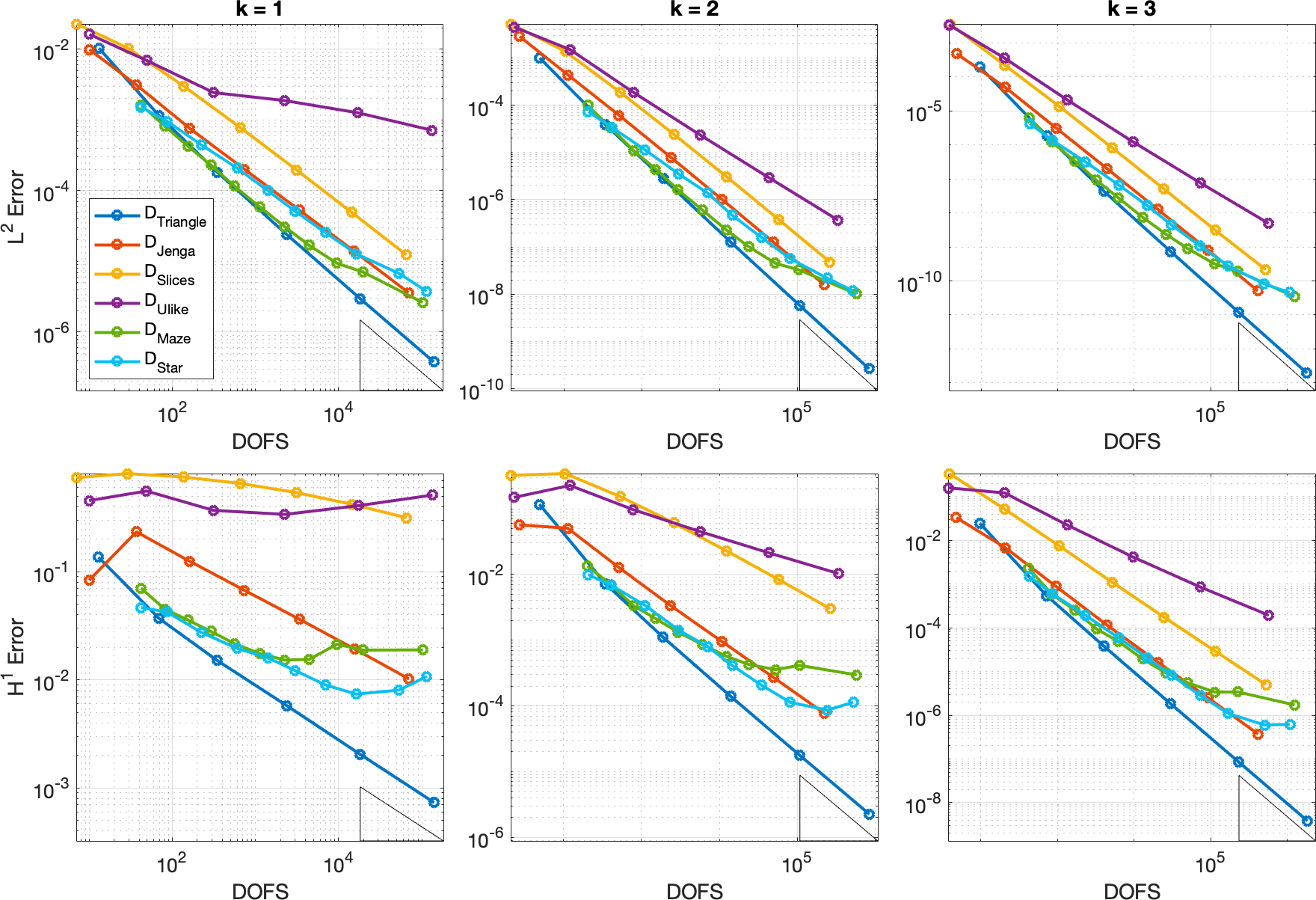} %% {Figures/performance1.png}
\caption{$L^2$-norm and $H^1$-seminorm of the approximation errors of the reference, hybrid and mirroring datasets for $k=1,2,3$. }
\label{fig:performance1}
\end{figure}

\begin{figure}[htbp]
\centering
	\includegraphics[width=\linewidth]{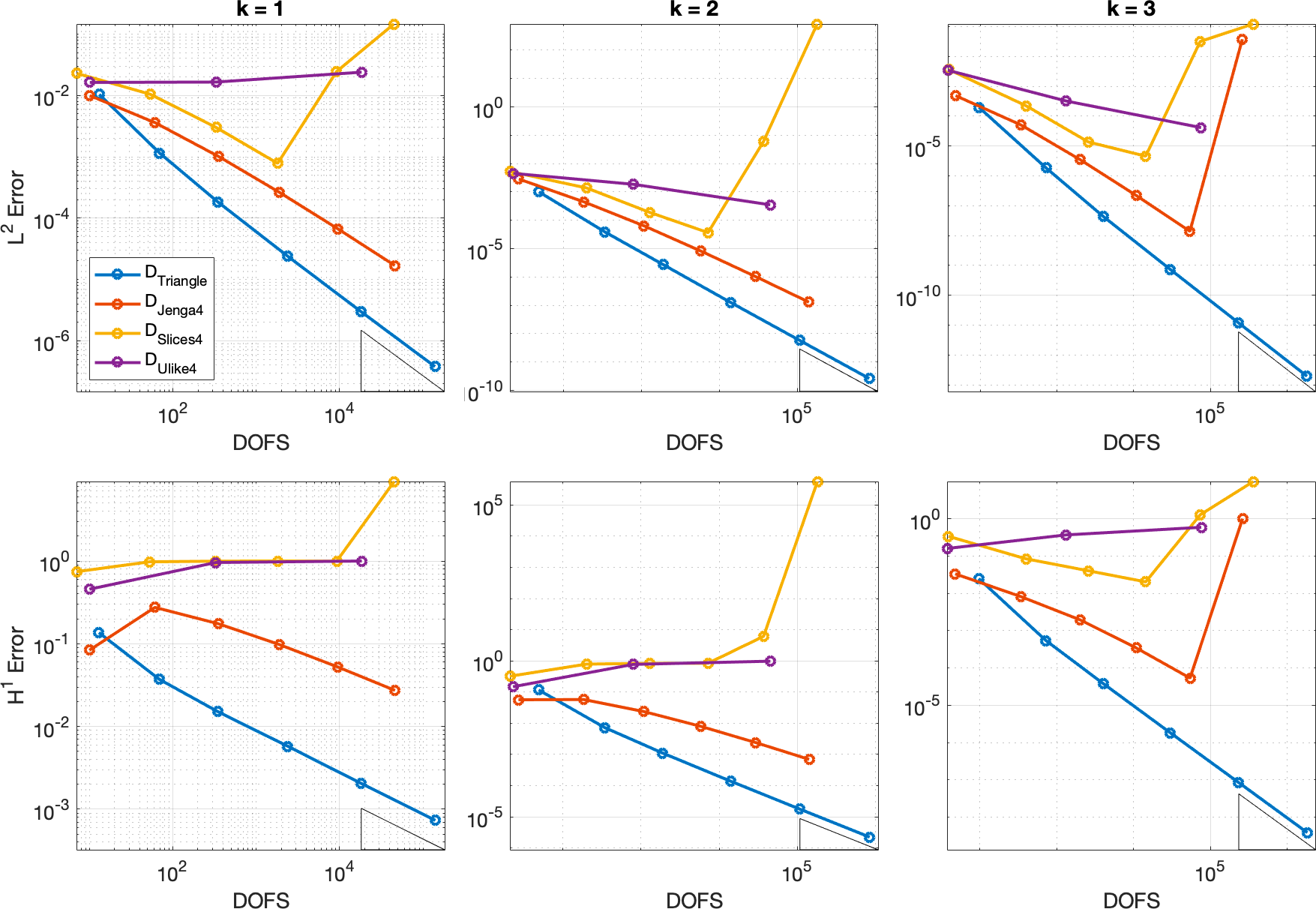} %% {Figures/performance4.png}
\caption{$L^2$-norm and $H^1$-seminorm of the approximation errors of the reference and multiple mirroring datasets for $k=1,2,3$. }
\label{fig:performance4}
\end{figure}

We solved the discrete Poisson problem \eqref{eq:Poisson:weak} with the VEM \eqref{eq:VEM} described in Section~\ref{sec:VEM} for $k=1,2,3$ over each mesh of each of the datasets defined in Section~\ref{subsec:datasets}, using as groundtruth the function 
\begin{equation}
    u(x,y)=\frac{\sin(\pi x)\sin(\pi y)}{2\pi^2}, \hspace{0.3cm} (x,y)\in\Omega=(0,1)^2.
\end{equation}
This function has homogeneous Dirichlet boundary conditions, and this choice was appositely made to prevent the boundary treatment from having an influence on the approximation error.
In Fig.~\ref{fig:performance1} and Fig.~\ref{fig:performance4} we plot the relative $L^2$-norm
$\norm{u-u_h}{0,\Omega}/\norm{u}{0,\Omega}$ and the relative $H^1$-seminorm $\snorm{u-u_h}{1,\Omega}/\snorm{u}{1,\Omega}$ (also called \textit{discrete energy norm}) of the approximation error as the number of DOFs increases (that is, as $n\to N$).\\

We also consider the condition numbers of matrices \textbf{G} and \textbf{H} (with the notation adopted in \cite{BeiraodaVeiga-Brezzi-Marini-Russo:2014}) as numerical indicators of the good behaviour of the method, and identities $|\matPin{k}\textbf{D}-\textbf{I}|=0$ and $|\matPiz{k}\textbf{D}-\textbf{I}|=0$ as an estimate of the approximation error produced by projectors $\Pin{k}$ and $\Piz{k}$, represented by matrices $\matPin{k}$ and $\matPiz{k}$, respectively. 
The computation of the projectors is obviously affected by the condition numbers of \textbf{G} and \textbf{H}, but the two indicators are not necessarily related.
All of these quantities are computed element-wise and the maximum value among all elements of the mesh is selected.
Condition numbers and identity values for $k=1,2,3$ are reported in Table~\ref{table:numerical_performance} (for $k<3$ we have $\matPiz{k}=\matPin{k}$).\\

First, the reference dataset $\Dtriangle$ guarantees for the correctness of the VEM, as it performs perfectly according to the theoretical results both in $L^2$ and in $H^1$ norms (the slopes being indicated by the triangles) for all $k$ values, maintaining reasonable condition numbers and optimal errors on the projectors $\Piz{k}$ and $\Pin{k}$.\\

For the hybrid datasets $\Dstar$ and $\Dmaze$, errors decrease at the correct rate for most of the meshes, and only start deflecting for very high numbers of DOFs and very complicated meshes.
These deflections are not due to numerical problems, as in both datasets we have cond(\textbf{G}) $<10^6$ and cond(\textbf{H}) $<10^9$, which are still reasonable values.
Projectors seem to work properly: $|\matPin{k}\textbf{D}-\textbf{I}|$ remains below $10^{-8}$ and $|\matPiz{k}\textbf{D}-\textbf{I}|$ below $10^{-7}$.
In a preliminary stage of this work, we obtained similar plots (not reported here) using other hybrid datasets built in the same way, with polygons surrounded by triangles.
In particular, we did not see big differences when starting with the other initial polygons of Benchmark \cite{attene2019benchmark}, cf. the construction discussed in ``Hybrid datasets'' in Section~\ref{subsec:datasets}.\\

On the meshes from ``Mirroring datasets'', $A_n$ or $e_n$ may scale non-uniformly, as reported in Table~\ref{table:ratios} (indeed, they can scale exponentially).
This reflects to cond(\textbf{G}) and cond(\textbf{H}), which grow up to $10^{10}$ and $10^{14}$ for $\Djenga$ in the case $k=3$.
Nonetheless, the discrepancy of the projectors identities remains below $10^{-5}$, which is not far from what happened with $\Dmaze$ and $\Dstar$.
Dataset $\Djenga$ exhibits an almost perfect convergence rate, even though $L^2$ and $H^1$ errors are bigger in magnitude than the ones measured for hybrid datasets;
$\Dslices$ shows even bigger errors and a non-optimal convergence rate, and $\Dulike$ is the dataset with the poorest performance, but still converges at a decent rate for $k>1$.\\

In the setting of ``Multiple mirroring datasets'', all datasets diverge badly (see Fig.~\ref{fig:performance4}), and this is principally due to very poor conditioning in the matrices involved in the calculations (see Table~\ref{table:numerical_performance}).
Dataset $\DjengaM$ and $\DslicesM$ maintain a similar trend to the ones in Fig.~\ref{fig:performance1} until numerical problems cause cond(\textbf{G}) and cond(\textbf{H}) to explode up to over $10^{30}$ for $\DjengaM$ and $10^{18}$ for $\DslicesM$.
In these conditions, projection matrices $\matPin{k}$ and $\matPiz{k}$ become meaningless and the method diverges.
The situation slightly improves for $\DulikeM$: cond(\textbf{H}) is still $10^{16}$, but the discrepancy of $\matPin{k}$ and $\matPiz{k}$ remain acceptable.
As a result, $\DulikeM$ does not properly explode, but the approximation error and the convergence rate are much worse than those seen in Fig.~\ref{fig:performance1}.\\

\begin{table}[htbp]
\centering
\caption{Summary of numerical performance for all datasets. We report the $\log_{10}$ of the original values for the condition number of \textbf{G} and \textbf{H} and the discrepancy of projection matrices $\matPin{k}$ and $\matPiz{k}$. Note that for $k<3$ we have $\matPiz{k}=\matPin{k}$.}
\label{table:numerical_performance}
%\resizebox{\columnwidth}{!}{
\begin{tabular}{c|ccc|ccc|ccc|ccc|ccc|ccc|ccc|ccc|ccc}

\hline
\multicolumn{1}{c}{\textbf{dataset}} &
\multicolumn{3}{c}{$\Dtriangle$} & \multicolumn{3}{c}{$\Dmaze$} &  \multicolumn{3}{c}{$\Dstar$} & \multicolumn{3}{c}{$\Djenga$} &  \multicolumn{3}{c}{$\Dslices$} & \multicolumn{3}{c}{$\Dulike$} & \multicolumn{3}{c}{$\DjengaM$} &  \multicolumn{3}{c}{$\DslicesM$} &  \multicolumn{3}{c}{$\DulikeM$}\\

\hline\noalign{\smallskip}
$k$ & 1 & 2 & 3 & 1 & 2 & 3 & 1 & 2 & 3 & 1 & 2 & 3 & 1 & 2 & 3 & 1 & 2 & 3 & 1 & 2 & 3 & 1 & 2 & 3 & 1 & 2 & 3\\

\noalign{\smallskip}\hline\noalign{\smallskip}
cond(\textbf{G}) 
& 0 & 2 & 5 % Triangle
& 2 & 3 & 6 % Maze
& 1 & 3 & 6 % Star
& 1 & 5 & 10 % Jenga
& 2 & 4 & 6 % Slices
& 1 & 4 & 7 % Ulike
& 6 & 18 & 31 % Jenga4
& 6 & 8 & 10 % Slices 4
& 2 & 6 & 11\\ % Ulike4

cond(\textbf{H})
& 2 & 5 & 7 % Triangle
& 2 & 5 & 8 % Maze
& 3 & 6 & 9 % Star
& 4 & 9 & 14 % Jenga
& 2 & 8 & 10 % Slices
& 3 & 7 & 10 % Ulike
& 13 & 26 & 39 % Jenga4
& 2 & 15 & 18 % Slices 4
& 5 & 10 & 16\\ % Ulike4

$|\matPin{k}\textbf{D}-\textbf{I}|$
& -13 & -11 & -9 % Triangle
& -12 & -10 & -8 % Maze
& -12 & -10 & -8 % Star
& -12 & -8 & -5 % Jenga
& -12 & -10 & -9 % Slices
& -13 & -10 & -8 % Ulike
& -9 & 3 & 13 % Jenga4
& -8 & -6 & -5 % Slices 4
& -13 & -8 & -5\\ % Ulike4  

$|\matPiz{k}\textbf{D}-\textbf{I}|$
&  &  & -10 % Triangle
&  &  & -8 % Maze
&  &  & -7 % Star
&  &  & -5 % Jenga
&  &  & -5 % Slices
&  &  & -7 % Ulike
&  &  & 20 % Jenga4
&  &  & 8 % Slices 4
&  &  & -4\\ % Ulike4

\noalign{\smallskip}\hline
\end{tabular}
%}
\end{table}

As a preliminary conclusion, by simply looking at the previous plots we observe that the relationship is not particularly strong between the geometrical assumptions respected by a certain dataset and the performance of the VEM.
In fact, we obtained reasonable results with meshes violating several assumptions.
\section{Mesh Quality Indicator}
\label{sec:indicator}
We now aim at defining a mesh quality indicator, that is, a scalar function capable of providing insights on the behaviour of the VEM over a particular sequence of meshes, before actually computing the approximated solutions.

\subsection{Definition}
\label{subsec:indicator:definition}
\begin{figure}[htbp]
\centering
	\includegraphics[width=.3\linewidth]{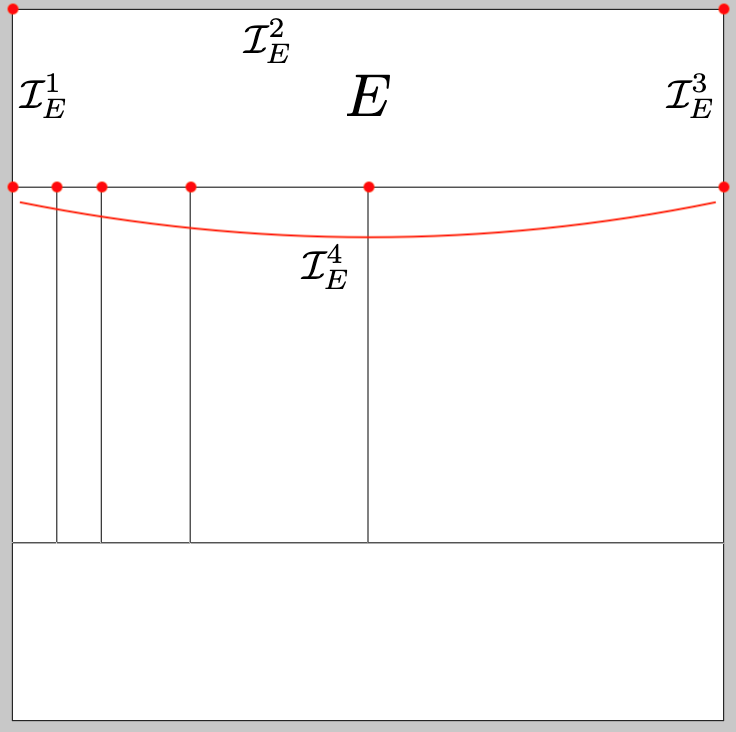} %% {Figures/assumption-g4-jenga.png}
\caption{One-dimensional mesh $\calI_\P=\{\calI_\P^1, \calI_\P^2, \calI_\P^3, \calI_\P^4\}$ for the top bar $\P$ of a $\Djenga$ base mesh.}
\label{fig:rho4}
\end{figure}

We start from the geometrical assumptions defined in Section~\ref{sec:geometrical:assumptions}.
Even if we proved them not to be strictly necessary for the convergence of the method, they can still be good indicators for the general quality of a sequence of meshes.
From each geometrical assumption \textbf{Gi}, \textbf{i }$=1,\ldots, 4$, we derived a scalar function $\varrho_i: \{\P \subset \Omega_h\} \to [0,1]$ defined element-wise, which measures how well a polygon $\P\in \Omega_h$ meets the requirements of \textbf{Gi} from 0 ($\P$ does not respect \textbf{Gi}) to 1 ($\P$ fully respects \textbf{Gi}).\\
\begin{align}
\varrho_1(\P) &= \frac{k(\P)}{|\P|} = 
    \begin{cases} 
    1 & \mbox{if $\P$ is convex} \\
    \in (0,1) & \mbox{if $\P$ is concave and star-shaped} \\
    0 & \mbox{if $\P$ is not star-shaped} \\
    \end{cases}\\
\varrho_2(\P) &= \frac{\min(\sqrt{|\P|}, \ \min_{e \in \partial \P}|e|)}{\max(\sqrt{|\P|}, \ h_\P)} \\
\varrho_3(\P) &= \frac{3}{\#\left\{ e \in \partial \P \right\}} = 
    \begin{cases} 
    1 & \mbox{if $\P$ is a triangle} \\
    \in (0,1) & \mbox{otherwise} \\
    \end{cases}\\
\varrho_4(\P) &= \min_i \frac{\min_{e\in \calI_\P^i} |e|}{\max_{e\in \calI_\P^i} |e|} \ \ %\mbox{where $\calI_\P = \cup_i \calI_\P^i$ is the one dimensional mesh related to $\partial \P$}.
\end{align}
The operator $k(\P)$ in $\varrho_1$ measures the area of the \textit{kernel} of a polygon $\P$, defined as the set of points in $\P$ from which the whole polygon is visible.
Therefore, $\varrho_1(\P)$ can be interpreted as an estimate of the value of the constant $\rho$ from assumption~\textbf{G1} on the polygon $\P$.
Similarly, the function $\varrho_2$ returns an estimate of the constant $\rho$ introduced in \textbf{G2}, expressed trough the ratio $|\E|/h_\P$, with the insertion of the quantity $\sqrt{|\P|}$ in order to avoid pathological situations.
Function $\varrho_3$ is a simple counter of the number of edges of a polygon, which penalizes elements with numerous edges as required by \textbf{G3}.
Last, we recall from Section~\ref{sec:geometrical:assumptions} that the boundary of a polygon $\P$ can be considered as a one-dimensional mesh $\calI_\P$, which can be subdivided into a number of disjoint sub-meshes $\calI_\P^1, \ldots, \calI_\P^N$, each one containing possibly more than one edge of $\P$.
In practice, we consider as a sub-mesh the collection of all edges whose vertices lie on the same line.
For example, as shown in Fig.~\ref{fig:rho4}, the boundary of the top bar $\P$ in the base mesh of $\Djenga$ is represented by a mesh $\calI_\P=\{\calI_\P^1, \calI_\P^2, \calI_\P^3, \calI_\P^4\}$, where $\calI_\P^1, \calI_\P^2$ and $\calI_\P^3$ contain, respectively, the left, top and right edge of $\P$, while $\calI_\P^4$ contains all the aligned edges in the bottom of $\P$.
Function $\varrho_4$ returns the minimum ratio between the smallest and the largest element in every $\calI_\P$, that is a measure of the quasi-uniformity of $\calI_\P$ imposed by \textbf{G4}.

\medskip
Combining together $\varrho_1, \varrho_2$, $\varrho_3$ and $\varrho_4$, we define a global function $\varrho:\{\Th\}_h\to [0,1]$ which measures the overall quality of a mesh $\Omega_h$.
Given a dataset $\calD$, we can study the behaviour of $\varrho(\Th)$ for $\Th\in\calD$ and determine the quality of the dataset through the refinement process.
In particular, we chose the formula $\varrho_1\varrho_2 + \varrho_1\varrho_3 + \varrho_1\varrho_4$ as it reflects the way in which the relative assumptions are typically imposed: 
\textbf{G1} and \textbf{G2}, \textbf{G1} and \textbf{G3} or \textbf{G1} and \textbf{G4} (but not, for instance, \textbf{G2} and \textbf{G3} simultaneously):

\begin{equation}
\varrho(\Omega_h) = \sqrt{\frac{1}{\#\left\{ \P \in \Omega_h \right\}} \ \sum_{\P\in\Omega} \frac{\varrho_1(\P)\varrho_2(\P) + \varrho_1(\P)\varrho_3(\P)+ \varrho_1(\P)\varrho_4(\P)}{3}}.
\label{eq:indicator}
\end{equation}
We have $\varrho(\Omega_h)=1$ if and only if $\Omega_h$ is made only of equilateral triangles, $\varrho(\Omega_h)=0$ if and only if $\Omega_h$ is made only of non star-shaped polygons, and $0<\varrho(\Omega_h)<1$ otherwise.
All indicators $\varrho_1, \varrho_2$, $\varrho_3$ and $\varrho_4$, and consequently $\varrho$, only depend on the geometrical properties of the mesh elements;
therefore their values can be computed before applying the VEM, or any other numerical scheme.\\
We point out that this approach is easily upgradeable to future developments:
whenever new assumptions on the features of a mesh should come up, one simply needs to introduce in our framework a new function $\varrho_i$ 
that measures the violation of the new assumption and insert it into the formulation of the general indicator $\varrho$ in equation~(\ref{eq:indicator}).
\subsection{Results}
\label{subsec:indicator:results}
\begin{figure}[htbp]
\centering
\begin{tabular}{c c}
	\includegraphics[width=.4\linewidth]{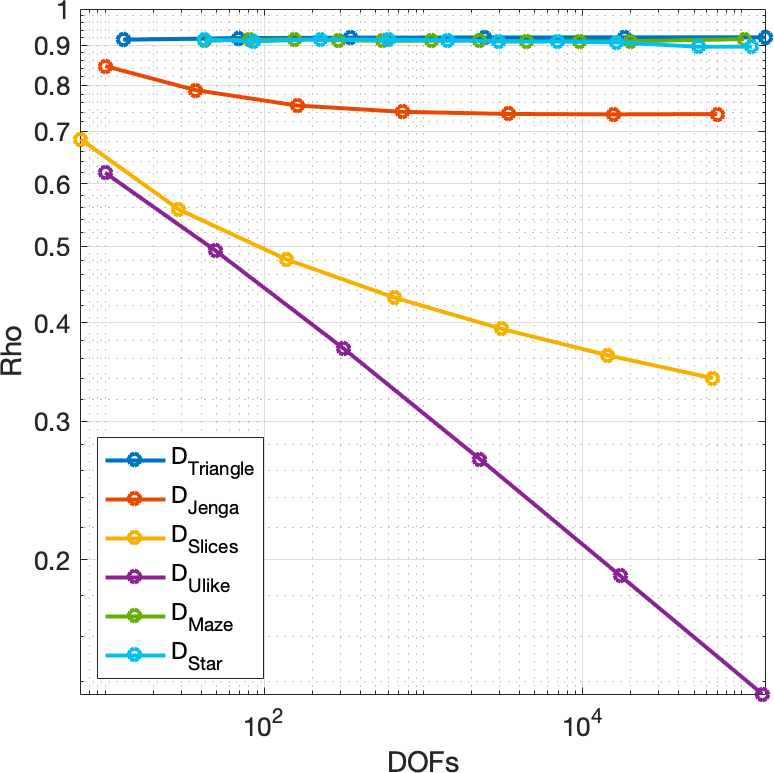} %% {Figures/indicator1.png} & \qquad
	\includegraphics[width=.4\linewidth]{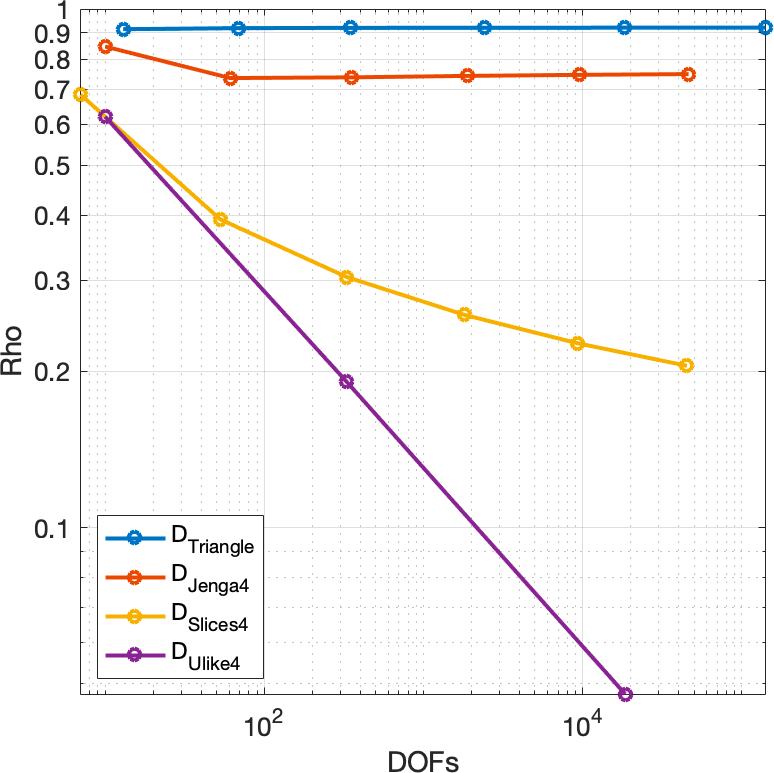} %% {Figures/indicator4.png} \\
    (a) & (b)\\
\end{tabular}
\caption{Indicator $\varrho$ for all datasets.}
\label{fig:indicator}
\end{figure}

We evaluated the indicator $\varrho$ over the datasets defined for this work; results are shown in Fig.~\ref{fig:indicator}.
If we compare Fig.~\ref{fig:indicator}(a) and \ref{fig:indicator}(b) with Fig.~\ref{fig:performance1} and \ref{fig:performance4} respectively, we can look for a correspondence between the behaviour of $\varrho$ on a dataset, computed before solving the problem, and the approximation error actually produced by that dataset.
Clearly, as $\varrho$ does not depend on the polynomial degree $k$ nor on the type of norm used, we will compare it to an average of the plots for the different $k$ values and for the different norms ($L^2$ and $H^1$).

\medskip
We preliminarily observe that, for an ideal dataset made by meshes containing only equilateral triangles, $\varrho$ would be constantly equal to 1.
We assume this value as a reference for the other datasets: 
the closer is $\varrho$ on a dataset to the line $y=1$, the smaller is the approximation error that we expect that dataset to produce.
Similarly, the more negative is the $\varrho$ slope, the worse is the convergence rate that we expect over that dataset.

For meshes belonging to $\Dtriangle$, $\varrho$ is almost constant and very close to 1, thus foreseeing the excellent convergence rates and the low errors seen in every sub-figure of Fig.~\ref{fig:performance1}.
The plots for $\Dmaze$ and $\Dstar$ in Fig.~\ref{fig:indicator}(a) are close to $\Dtriangle$, hence we expect them to behave similarly.
This is confirmed by Fig.~\ref{fig:performance1}: $\Dmaze$ and $\Dstar$ are almost coincident and very close to $\Dtriangle$ until the very last meshes, especially in the $L^2$ plots.

The $\Djenga$ plot in Fig.~\ref{fig:indicator}(a) anticipates a perfect convergence rate but greater error values with respect to the previous three, and again this behaviour is respected in Fig.~\ref{fig:performance1}.
The curve relative to $\Dslices$ in Fig.~\ref{fig:indicator}(a) is quite distant from the ideal value of 1.
Importantly, it keeps decreasing from mesh to mesh, even if the plot allows us to assume that it may flatten within a couple more meshes.
Looking at Fig.~\ref{fig:performance1}, we notice that this dataset produces an error significantly higher than the previous ones ($\Dtriangle, \Dmaze, \Dstar, \Djenga$), and in some cases the $H^1$ error convergence rate is significantly lower than the theoretical estimate.
Last, the $\varrho$ values in Fig.~\ref{fig:indicator}(a) predict huge errors and a completely wrong convergence rate for $\Dulike$.
This dataset is actually the one with the worst performance in Fig.~\ref{fig:performance1}, where it does not even always converge (see the case $k=1$, $H^1$ seminorm).

As far as multiply refined datasets are concerned, we notice that, since it only depends on the geometry of the elements, $\varrho$ is not affected by numerical errors.
The $\varrho$ plot for $\DjengaM$ in Fig.~\ref{fig:indicator}(b) is very similar to the plot obtained for $\Djenga$ in \ref{fig:indicator}(a), therefore we should expect $\DjengaM$ in Fig.~\ref{fig:performance4} to perform similarly to $\Djenga$ in Fig.~\ref{fig:performance1}.
This is actually the case at least until the last mesh for $k=3$, when numerical problems appear which $\varrho$ is not able to predict.
Also dataset $\DslicesM$ has a similar trend to $\Dslices$ but decreases faster, reaching a $\varrho$ value of $\sim0.2$ instead of $\sim0.34$ within a smaller number of meshes.
As above, $\DslicesM$ performs similarly to $\Dslices$ until condition numbers explode, in the last two meshes for every value of $k$.
Last, the $\varrho$ plot of $\DulikeM$ is significantly worse than the one of $\Dulike$ (and than any other), both in terms of distance from $y=1$ and slope.
In Fig.~\ref{fig:performance4} we can observe how, even if $\DulikeM$ does not properly explode (as it suffers less from numerical problems, cf. Table~\ref{table:numerical_performance}), the approximation error and the convergence rate are the worse among all the considered datasets.

\medskip
Summing up these results, we %can 
conclude that indicator $\varrho$ is able, up to a certain accuracy, to predict the behaviour of the VEM over the considered datasets, both in terms of error magnitude and convergence rate.
The prediction may be inaccurate in presence of very similar performance (the case of $\Dmaze$ and $\Dstar$), or in extreme situations in which the numerical problems become so significant to overcome any influence that the geometrical features of the mesh could have on the performance (the last meshes of $\DjengaM$ and $\DslicesM$).
\section{Conclusions}
\label{sec:conclusions}
In this work, we collected the regularity assumptions that are used in the literature to guarantee the convergence and the error estimates in the $\LTWO$ and $\HONE$ norms for the VEM.
These conditions allow a great flexibility for the type and variety of polygons to be used in a mesh, but they still seem overestimated.
Experimentally, we verified that the VEM works, with a good convergence rate, also on meshes and datasets that strongly violate these assumptions.
We also built examples of datasets for which, violating significantly the regularity assumptions, the VEM shows a convergence rate suboptimal or diverges.
Finally, we introduced new indicators to represent how much the regularity hypothesis are violated by a tessellation and combined these indicators in a single score, aimed at estimating how a dataset can be expected to be performing in the solution of the VEM.
The results obtained are encouraging, showing a satisfactory correlation between the errors and this indicator.
Consequently, our approach provides an experimental score that is able to predict if a tessellation of a domain can be critical for the VEM.

As possible future developments, we are interested in refining the regularity indicator here proposed,
for example, to deduce new decomposition rules of a domain with possible applications to mesh
generators, or to adaptive coarsening/refinement algorithms.
We are also experimenting similar indicators to evaluate the properties of polyhedral meshes.

\appendix

\section{Main convergence results in the literature of the VEM}
\label{app:error}
This appendix is a short overview of the main results on convergence analysis from the VEM
literature.
For each selected paper, we report (where available) results for
abstract energy error, $\HONE$ error estimate and $\LTWO$ error
estimate, highlighting the geometrical assumptions considered.
We may have changed the notation in a few points and introduced some
very minor modifications in the theorem statements for consistency
with our paper.

\medskip
%%%%%%%%%%%%%%%%%%%%%%%%%%%
\subsection{``Basic Principles of Virtual Elements Methods''
  \cite{BeiraodaVeiga-Brezzi-Cangiani-Manzini-Marini-Russo:2013}}
\label{app:error:basic}
\noindent
This work is the very first paper about the VEM, where this method
was introduced.
The original formulation adopts the regular conforming virtual element
space, which we still denote by $\Vhk$ as in \eqref{eq:Vhk:def} and \eqref{eq:VhkP:def} with a small abuse of notation:
\begin{subequations}
  \label{eq:Vhk:regular:conforming}
  \begin{align}
    \Vhk:= 
    &\{\vsh\in\HONE(\Omega):\ \restrict{\vsh}{\P}\in\Vhk(\P) \ \forall\P\in\Th\},
    \intertext{where}
    \Vhk(\P):= 
    &\{\vsh\in\HONE(\P):\ 
    \restrict{\vsh}{\partial\P}\in C^0(\partial\P), \ 
    \restrict{\vsh}{\E}\in\PS{k}(\E) \ \forall\E\in\partial\P, \ 
    \Delta\vsh\in\PS{k-2}(\P)\},
  \end{align}
\end{subequations}
and the \textit{dofi-dofi} formulation $\SPhdd$ defined in \eqref{eq:stab:dofi-dofi} is introduced for the stabilization bilinear form.\\
Although not explicitly used to derive the following abstract result for the convergence in the energy norm, the Authors
introduce the mesh regularity assumptions \textbf{G\ref{g1}} and
\textbf{G\ref{g2}} and the concept of \emph{simple polygon}, which is
a connected polygonal element with a nonintersecting boundary made of
straight edges.
This setting is the general and widely adopted framework of the
virtual element formulation in many successive papers.
Moreover, a broken $\HONE$-seminorm is introduced, for functions $\vs\in\HONE(\Th)$:
\begin{align}
    \snorm{\vs}{\hh,1}:=\left(\sum_{\P\in\Th} \snorm{\nabla \vs}{0,\P}^2\right)^{1/2}.
    \label{eq:h1-broken}
\end{align}
\begin{theo}[abstract energy error]
\label{theo:basic:abstract}
  Under the \textit{k-consistency} and \textit{stability} assumptions
  defined in Section~\ref{sec:VEM}, cf. \eqref{eq:k-consistency} and \eqref{eq:stability}, the discrete
  problem has a unique solution $\ush$.
  Moreover, for every approximation $\usI\in \Vhk$ of $\us$ and every
  approximation $\us_\pi$ of $\us$ that is piecewise in $\PS{k}(\Th)$,
  we have
  \begin{equation}
    \snorm{\us-\ush}{1,\Omega} \leq
    \Cs(\snorm{\us-\usI}{1,\Omega} +
    \snorm{\us-\us_{\pi}}{\hh,1} +
    \mathfrak{F}_{\hh}),
  \end{equation}
  where $\Cs$ is a constant depending only on $\alpha_*$ and
  $\alpha^*$ (the constants in \eqref{eq:stability}), and, for any $\hh$,
  $\mathfrak{F}_{\hh}=\snorm{\fs-\fsh}{{\Vhk}^\prime}$ is the
  smallest constant such that
  \begin{equation*}
    (\fs,\vs)-\langle\fsh-\fs,\vs\rangle \leq
    \mathfrak{F}_{\hh}\snorm{\fs}{1}
    \qquad\forall\vs\in\Vhk.
  \end{equation*}
\end{theo}
The Authors claim that an $\LTWO$ error estimate of the
convergence rate can be derived with the usual duality argument
techniques.

\medskip
%%%%%%%%%%%%%%%%%%%%%%%%%%%%%%%%%%%
\subsection{``Equivalent projectors for virtual element methods'' \cite{Ahmad-Alsaedi-Brezzi-Marini-Russo:2013}}
\label{app:error:equivalent}
\noindent
In this paper $\Vhk(\P)$ is replaced by the enhanced VEM space \eqref{eq:VhkP:def} adopted in our work (in the paper it is called ``modified VEM space'') and the \textit{dofi-dofi} stabilization is adopted.
Under the geometrical assumptions \textbf{G\ref{g1}} and \textbf{G\ref{g2}}, $\HONE$ and $\LTWO$ error estimates are provided; while for the abstract energy error, Theorem \ref{theo:basic:abstract} is reported.
\begin{theo}[$\HONE$ error estimate]
  Assuming~\textbf{G\ref{g1}},~\textbf{G\ref{g2}}, let the right-hand side
  $\fs$ belong to $H^{k-1}(\Omega)$, and that the exact solution $\us$ belong
  to $H^{k+1}(\Omega)$. Then
  \begin{equation}
    \norm{\us-\ush}{1,\Omega} \le \Cs|\hh|^k \snorm{\us}{k+1,\Omega}
  \end{equation}
  with $\Cs$ a positive constant independent of $\hh$.
\end{theo}
\begin{theo}[$\LTWO$ error estimate]
  Assuming~\textbf{G\ref{g1}},~\textbf{G\ref{g2}} and with $\Omega$ convex,
  let the right-hand side $\fs$ belong to $H^k(\Omega)$, and that the exact
  solution $\us$ belong to $H^{k+1}(\Omega)$. Then
  \begin{equation}
    \norm{\us-\ush}{0,\Omega} + 
    |\hh|\norm{\us-\ush}{1,\Omega} \le
    \Cs|\hh|^{k+1}\snorm{u}{k+1,\Omega},
  \end{equation}
  with $\Cs$ a constant independent of $\hh$.
\end{theo}

\medskip
%%%%%%%%%%%%%%%%%%%%%%%%%%%%%%%%%%%
\subsection{``Stability analysis for the virtual element method'' \cite{Beirao-Lovadina-Russo:2016}}
\label{app:error:stability}
\noindent
This paper is based on the regular conforming VEM space \eqref{eq:Vhk:regular:conforming} defined in \cite{BeiraodaVeiga-Brezzi-Cangiani-Manzini-Marini-Russo:2013}.
A new abstract energy error estimate is deduced, and the $\HONE$ error is studied considering two different stabilization techniques.
The Authors also introduce new analytical assumptions on the bilinear form $\ash(\cdot,\cdot)$, replacing \eqref{eq:stability}:
\begin{subequations}
\begin{align}
    \ashP(\vsh,\vsh) 
    &\le \Cs_1(\P)\Tnorm{\vsh}{\P}^2, 
    \text{ for all } \vsh\in\Vhk(\P);\\
    \Tnorm{\qs}{\P}^2 
    &\le \Cs_2(\P)\ashP(\qs,\qs), 
    \text{ for all } \qs\in\PS{k}(\P),
\end{align}
\label{eq:stability2}
\end{subequations}
being $\Tnorm{\cdot}{\P}$ a discrete semi-norm induced by the stability term and $\Cs_1(\P), \Cs_2(\P)$ positive constants which depend on the shape and possibly on the size of $\P$.
Differently than the standard analysis of \cite{BeiraodaVeiga-Brezzi-Cangiani-Manzini-Marini-Russo:2013} where a kind of bound \eqref{eq:stability2}(b) is assumed for every $\vsh\in\Vhk(\P)$, here the estimate is only required for the polynomials $\qs\in\PS{k}(\P)$. 
Thus, even when $\Cs_1(\P)$ and $\Cs_2(\P)$ can be chosen independent of $\P$, on $\Vhk(\P)$ the semi-norm induced by the stabilization term may be stronger than the energy $\ashP(\cdot,\cdot)^{1/2}$.\\
For the following theorem, from the constants in \eqref{eq:stability2} the Authors derive the quantities:
\begin{align*}
    \tilde{\Cs}(\hh) = \max_{\P\in\Th}\{1,\,\Cs_2(\P)\},\quad
    \Cs_1(\hh) = \max_{\P\in\Th}\{\Cs_1(\P)\},\quad
    \Cs^*(\hh) = \frac{1}{2}\max_{\P\in\Th}\{\min\{1,\Cs_2(\P)^{-1}\}\},
\end{align*}
\begin{theo}[abstract energy error]
  Under the \textit{stability} assumptions \eqref{eq:stability2},
  let the continuous solution $\us$ of the problem satisfy
  $\restrict{\us}{\P}\in \calV_\P$ for all $\P\in\Th$, where $\calV_\P \subseteq \Vhk(\P)$ is a subspace of sufficiently regular functions.
  Then, for every $\usI\in \Vhk$ and for every $\us_\pi$ such that
  $\us_{\pi|\P}\in\PS{k}(\P)$, the discrete solution $\ush$ satisfies
  \begin{equation}
    \snorm{\us-\ush}{1,\Omega} \lesssim
    \Cs_{\text{err}}(\hh)\ (\mathfrak{\Fs}_{\hh}
    + \Tnorm{\us-\usI}{}
    + \Tnorm{\us-\us_{\pi}}{}
    + \snorm{\us-\usI}{1,\Omega}
    + \snorm{\us-\us_{\pi}}{\hh,1}),
  \end{equation}
  where the constant $\Cs_{\text{err}}(\hh)$ is given by
  \begin{equation*}
    \Cs_{\text{err}}(\hh) = \max\left\{
    1,\,\tilde{\Cs}(\hh)\Cs_1(\hh),\tilde{\Cs}(\hh)^{3/2}\sqrt{\Cs^*(\hh)\Cs_1(\hh)}
    \right\}.
  \end{equation*}
\end{theo}
The Authors consider the stability term $\SPh(\cdot,\cdot)$ as the sum of two
contributions: the first, $\Ss_h^{\partial\P}$, involving the boundary
degrees of freedom; the second, $\Ss_h^{o\P}$, involving the internal
degrees of freedom.
It can be shown that, for the following results, we can restrict the
analysis to $\Ss_h^{\partial\P}$, which can be expressed in the
\textit{dofi-dofi} form $S_{\hh}^{\partial\P,\textrm{dd}}$ already defined in \eqref{eq:stab:dofi-dofi}, or in the \textit{trace} form proposed in \cite{Wriggers-Rust-Reddy:2016}:
\begin{align}
  S_{\hh}^{\partial\P,\textrm{tr}}(\vsh,\wsh) =
  \hh_{\P} \int_{\partial\P} \partial_s\vsh \partial_s\wsh ds,
  \label{eq:stab:trace}
\end{align}
where $\partial_s\vsh$ denotes the tangential derivative of $\vsh$ along $\partial\P$.
\begin{theo}[$\HONE$ error estimate with \textit{dofi-dofi} stabilization]
  Assuming~\textbf{G\ref{g1},~G\ref{g3}}, let $\us\in\HS{s}(\Omega)$, 
  $s>1$, be the solution of the problem with $\SPh=\Ss_h^{\partial\P,\textrm{dd}}$.
  Let $\ush$ be the solution of the discrete problem, then it holds
  \begin{align}
    \norm{\us-\ush}{1,\Omega} \lesssim \Cs(\hh)\hh^{s-1}\snorm{\us}{s,\Omega}\quad 1<\ss\le\ks+1,
  \end{align}
  with
  \[\Cs(\hh)=\max_{\P\in\Th}(\log(1+\hP/\hh_m(\P))),\]
  where $\hh_m(\P)$ denotes the length of the smallest edge of $\P$.
\end{theo}
\begin{corollary}
  Assuming~\textbf{G\ref{g1}} and \textbf{G\ref{g2}} instead, then
  $\cs(\hh)\lesssim 1$ and therefore
  \begin{equation*}
    \norm{\us-\ush}{1,\Omega}\lesssim\hh^{s-1}\snorm{\us}{s,\Omega}\qquad 1<s\le\ks+1.
  \end{equation*}
\end{corollary}        
\begin{theo}[$\HONE$ error estimate with \textit{trace} stabilization]
  Under Assumption~\textbf{G\ref{g1}}, let $\us\in\HS{s}(\Omega)$,
  $s>3/2$ be the solution of the problem with $\SPh=\Ss_h^{\partial\P,\textrm{tr}}$.
  %%%
  Let $\ush$ be the solution of the discrete problem, then it holds
  \begin{equation}
    \norm{\us-\ush}{1,\Omega}\lesssim\hh^{s-1}\snorm{\us}{s,\Omega}\quad 3/2<s\le k+1.
  \end{equation}
\end{theo}

\medskip
%%%%%%%%%%%%%%%%%%%%%%%%%%%%%%%%%%%
\subsection{``Some Estimates for Virtual Element Methods'' \cite{Brenner:2017:SEV}}
\label{app:error:estimates}
\noindent
In this paper, the enhanced VEM space is defined in a slightly different (but still equivalent) formulation from \eqref{eq:VhkP:def}:
\begin{equation}
    \label{eq:Vhk:enhanced2}
    \begin{split}
    \Vhk(\P) := \big\{ \vsh\in\HONE(\P):\
    &\restrict{\vsh}{\partial\P}\in\PS{k}(\partial\P),\\
    & \exists\,\qs_{\vsh}(=-\Delta\vsh)\in\PS{k}(\P)\text{~such~that~}
    \int_{\P}\nabla\vsh\cdot\nabla\wsh\dxv = \int_{\P}\qs_{\vsh}\wsh\dxv\quad\forall\wsh\in\HONEzr(\P), \\
    &\text{~and~}
    \Pi_k^{0,\P}\vsh - \Pi_k^{\nabla, \P}\vsh\in\PS{k-2}(\P) \big\}.
    \end{split}
\end{equation}
Different stabilization types are considered, but the convergence results in this case do not depend on the choice of $\SPh$.
The geometrical assumptions required throughout the article are \textbf{G\ref{g1}} and \textbf{G\ref{g2}}.
\begin{theo}[abstract energy error]
  Assuming~\textbf{G\ref{g1},~G\ref{g2}}, if
  $\fs\in\HS{s-1}(\Omega)$ for $1\le\ss\le\ks$, then there exists a
  positive constant $\Cs$ depending only on $k$ and $\rho$ from \textbf{G\ref{g1}} such that
  \begin{equation}
    \snorm{\us-\ush}{1,\Omega}\le\Cs(\inf_{\vs\in\Vhk}\snorm{\us-\vs}{1,\Omega}
    + \inf_{\ws\in\PS{k}(\Th)}\snorm{\us-\ws}{\hh,1}
    + \hh^s\snorm{\fs}{s-1,\Omega}).
  \end{equation}
\end{theo}
\begin{theo}[$\HONE$ error estimate]
  Assuming~\textbf{G\ref{g1},~G\ref{g2}}, if
  $\us\in\HS{s+1}(\Omega)$ for $1\le\ss\le\ks$, then there exists
  positive constants $\Cs_1$, $\Cs_2$ depending only on $k$ and $\rho$ from \textbf{G\ref{g1}} such that
  \begin{equation}
    \snorm{\us-\ush}{1,\Omega} + \snorm{\us-\Pin{k}\ush}{\hh,1}
    \leq\Cs_1\hh^s\snorm{\us}{s+1,\Omega}.
  \end{equation}
\end{theo} 
\begin{theo}[$\LTWO$ error estimate]
  Assuming~\textbf{G\ref{g1},~G\ref{g2}}, with $\Omega$
  convex, if $\us\in\HS{s+1}(\Omega)$ for for $1\le\ss\le\ks$, then
  there exists a positive constant $\Cs$ depending only on $\Omega$,
  $k$ and $\rho$ from \textbf{G\ref{g1}} such that
  \begin{equation}
    \norm{\us-\ush}{0,\Omega}\le\Cs\hh^{s+1}\snorm{\us}{s+1,\Omega}.
  \end{equation}
\end{theo}

\medskip
%%%%%%%%%%%%%%%%%%%%%%%%%%%%%%%%%%%
\subsection{``Virtual element methods on meshes with small edges or faces'' \cite{brenner2018virtual}}
\label{app:error:small}
\noindent
The Authors establish error estimates for virtual element methods on polygonal
or polyhedral meshes that can contain small edges $(d = 2)$ or small faces 
$(d = 3)$. The VEM space is the enhanced space formulated as in \eqref{eq:Vhk:enhanced2}, and the local stabilizing bilinear form is considered in the \textit{dofi-dofi} formulation $\SPhdd$ and in the \textit{trace} formulation $\SPhtr$ of \eqref{eq:stab:trace}. Also, the following constants are defined:
\begin{align}
  \label{eq:error:small}
  \calH := \sup_{\P\in\Th}\left(
  \frac{\max_{\E\in\partial\P}\hE}{\min_{\E\in\partial\P}\hE}
  \right), \qquad
  \alpha_{\hh}:= \begin{cases}
    \ln\left( 1+\calH \right) & \mbox{~with } \SPhdd\\
    1 & \mbox{~with } \SPhtr
  \end{cases}
\end{align}
The geometrical assumptions required throughout the article are \textbf{G\ref{g1}} and \textbf{G\ref{g3}}.
The Authors introduce a mesh-dependent energy norm $\norm{\cdot}{\hh} := \sqrt{\ash(\cdot,\cdot)}$ and a functional $\Xi_{\hh}:\Vhk\to\PS{k}(\Th)$ given by
\begin{align}
    \Xi_{\hh} = \begin{cases}
    \Piz{1} & \mbox{ if } k=1,2\\
    \Piz{k-1} & \mbox{ if } k\ge3.
    \end{cases}
\end{align}
\begin{theo}[abstract energy error]
  Assuming~\textbf{G\ref{g1},~G\ref{g3}}, let $\us$ and $\ush$ be the solutions of the continuous and discrete problems. We have:
  \begin{equation}
    \norm{\us-\ush}{\hh}\lesssim
    \inf_{\ws\in\Vhk} \norm{\us-\ws}{h} +
    \norm{\us-\Pin{k}\us}{\hh} +
    \sqrt{\alpha_{\hh}}\left(\norm{\us-\Pin{k}\us}{\hh,1} +
    \sup_{\ws\in\Vhk}\frac{(\fs,\ws-\Xi_{\hh}\ws)}{\snorm{w}{1,\Omega}}
    \right).
  \end{equation}
\end{theo}
\begin{theo}[$\HONE$ error estimate]
  Assuming~\textbf{G\ref{g1},~G\ref{g3}}, if the solution $\us$
  belongs to $\HS{s+1}(\Omega)$ for some $1\le\ss\le\ks$, we have:
  \begin{align}
    &\norm{\us-\ush}{\hh} \lesssim \sqrt{\alpha_{\hh}}\hh^s\snorm{u}{s+1,\Omega}, \ \textit{and}\\
    &\snorm{\us-\ush}{1,\Omega}
    + \sqrt{\alpha_{\hh}}\left[
      \snorm{\us-\Pin{k}\ush}{\hh,1}
      + \snorm{\us-\Piz{k}\us}{\hh,1}
      \right]
    \lesssim \alpha_{\hh}\hh^s\snorm{\us}{s+1,\Omega}.
  \end{align}
\end{theo}
\begin{theo}[$\LTWO$ error estimate]
  Assuming~\textbf{G\ref{g1},~G\ref{g3}}, if the solution $\us$ belongs to \\
  $\HS{s+1}(\Omega)$ for some $1\le\ss\leq\ks$, we have:
  \begin{equation}
    \norm{\us-\ush}{0,\Omega}\leq
    \Cs\ \alpha_{\hh}\hh^{s+1}
    \snorm{\us}{s+1,\Omega}.
  \end{equation}
\end{theo} 
The notation $A\lesssim B$ indicates that $A\le\Cs B$, with a positive constant $\Cs$ depending on the mesh regularity parameter $\rho$ of \textbf{G\ref{g1}} and the degree $k$ in the case of $\SPhtr$, and also on yhe maximum number of edges $N$ of \textbf{G\ref{g3}} in the case of $\SPhdd$.

\medskip
%%%%%%%%%%%%%%%%%%%%%%%%%%%%%%%%%%%
\subsection{``Sharper error estimates for Virtual Elements and a bubble-enriched version'' \cite{da2020sharper}}
\label{app:error:sharper}
\noindent
In this paper, it is shown that the $\HONE$ interpolation error $\snorm{\us-\usI}{1,\P}$
on each element $\P$ can be split into a boundary contribution and a
bulk contribution.
The idea is to decouple the polynomial order on the boundary and in
the bulk of the element.
Let $\ks_o$ and $\ks_{\partial}$ be two positive integers with
$\ks_o\ge\ks_{\partial}$ and let $\kv=(\ks_o,\ks_{\partial})$.
For any $\P\in\Th$ the Authors define the \textit{generalized virtual element
space}:
\begin{subequations}
\begin{align}
    V^h_{\kv} := 
    &\{\vsh\in\HONEzr(\Omega):
    \restrict{\vsh}{\P}\in) V^h_{\kv}(\P), \ \forall\P\in\Th\}, 
    \intertext{where}
    V^h_{\kv}(\P) := 
    &\{\vsh\in\HONEzr(\P):
    \restrict{\vsh}{\partial\P}\in\CS{0}(\partial\P), \
    \restrict{\vsh}{\E}\in\PS{\ks_\partial}(\E) \ \forall\E\in\partial\P, \
    \Delta\vsh\in\PS{\ks_o-2}(\P)\}.
\end{align}
\end{subequations}
For $\ks_o=\ks_{\partial}$, the space $V^h_{\kv}(\P)$ coincides with the regular virtual element space in \eqref{eq:Vhk:regular:conforming}.
Moreover, given a function $v \in \HONEzr \cap H^s(\Th)$, on each element $\P \in \Th$ the Authors define the interpolant function $\calI_hv$ as the solution of the following elliptic problem:
\begin{align*}
    \begin{cases} \Delta\calI_hv = \Pi^{0,\P}_{k_o-2} \Delta v& \mbox{in } \P \\
    \calI_hv = v_b & \mbox{on } \partial \P, \end{cases}
\end{align*}
where $v_b$ is the standard 1D piecewise polynomial interpolation of $v|_{\partial \P}$.
\begin{theo}[abstract energy error]
  Under Assumption~\textbf{G\ref{g1}}, let $u \in \HONEzr(\Th) \cap H^s(\Th)$ with $s > 1$ be the solution of the continuous problem and $\ush\in V^h_{\kv}$ be the solution of the discrete problem. Consider the functions
  \begin{equation*}
      e_h=\ush-\calI_hu, \ e_{\calI}=u-\calI_hu, \ e_{\pi}=u-u_{\pi}, \ e_u=u_{\pi}-\calI_hu,
  \end{equation*}
  where $u_{\pi} \in \PS{k_o}(\Th)$ is the piecewise polynomial approximation of $\us$ defined in Bramble-Hilbert Lemma. Then it holds that
  \begin{equation}
      \snorm{u-\ush}{1,\Omega}^2 + \alpha \ a_h(e_h,e_h) \lesssim \alpha^2 \sum_{\P\in \Th} h^2_\P \norm{f-f_h}{0,\P}^2 + \alpha^2 \snorm{e_\pi}{1,\Th}^2 + \alpha \snorm{e_\calI}{1,\Omega}^2 + \alpha \sum_{\P \in \Th}\sigma^\P
  \end{equation}
  where $\alpha$ is the coercivity constant and $\sigma^\P := \SPh ((I-\Pi^{\nabla,\P}_{k_0})e_u, (I-\Pi^{\nabla,\P}_{k_0})e_u)$.
\end{theo}
\begin{theo}[$\HONE$ error estimate with \textit{dofi-dofi} stabilization] 
  Assuming~\textbf{G\ref{g1},~G\ref{g4}}, let
  $\us\in \HONEzr(\Th)$ be the solution of the continuous problem and
  $\ush\in V^h_{\kv}$ be the solution of the discrete problem obtained
  with the \textit{dofi-dofi} stabilization.
  Assume moreover that $\us\in\HS{\bar{k}}(\Th)$ with
  $\bar{\ks}=\max\{\ks_o+1,\,\ks_{\partial}+2\}$ and
  $\fs\in\HS{\ks_o-1}$.
  Then it holds that
  \begin{equation}
    \snorm{\us-\ush}{1,\Omega}^2
    \lesssim\alpha\sum_{\P\in\Th}\left(
    (\alpha + \calN_{\P})^{1/2}\hP^{k_o}+\hh_{\partial\P}^{\ks_\partial}
    \right)^2,
  \end{equation}
  where $\hh_{\partial\P}$ denotes the maximum edge length,
  $\alpha$ is the constant defined in \eqref{eq:error:small}, and
  $\calN_{\P}$ is the number of edges in $\P$.
\end{theo}
\begin{theo}[$\HONE$ error estimate with \textit{trace} stabilization] 
  Under Assumption~\textbf{G\ref{g1}}, let $\us\in \HONEzr(\Th)$ be the
  solution of the continuous problem and $\ush\in V^h_{\kv}$ be the
  solution of the discrete problem obtained with the \textit{trace}
  stabilization.
  Assume moreover that $\us\in\HS{\bar{\ks}}(\Th)$ with
  $\bar{\ks}=\max\{\ks_o+1,\,\ks_{\partial}+2\}$ and
  $\fs\in\HS{\ks_o-1}$.
  Then it holds that
  \begin{equation}
    \snorm{\us-\ush}{1,\Omega}^2
    \lesssim\sum_{\P\in\Th}\left( \hP^{\ks_o} + \hh_{\partial\P}^{\ks_{\partial}}\right)^2.
  \end{equation}
\end{theo}
\section{Dataset generation}
\label{app:dataset}
In this appendix, we take a closer look at how the datasets presented in Section \ref{subsec:datasets} are built. 
All algorithms have been written using CinoLib \cite{livesu2019cinolib}.
We recall that a dataset is a finite mesh sequence $\calD=\{\Omega_n\}_{n=0,\ldots,N}$, ordered decreasingly with respect to the mesh size.
We also recall the definition of the quantities:
\begin{align*}
A_n=\frac{\max_{\P\in \Omega_n} |\P|}{\min_{\P\in \Omega_n} |\P|}
\quad\textrm{and}\quad 
e_n=\frac{\max_{e\in \Omega_n} |e|}{\min_{e\in \Omega_n} |e|}, 
\quad \textrm{for } n=0,\ldots,N.
\end{align*}

%%%%%%%%%%%%%%%%%%%%%%%%%%%%%%%%%%%%%
\subsection{Reference dataset}
\label{app:dataset:reference}
The first dataset, $\Dtriangle$, contains only triangular meshes that are built by inserting a number of points in the domain, and connecting them in a Delaunay triangulation.
The point set is defined through the \textit{Poisson Disk Sampling} algorithm proposed in \cite{bridson2007fast}, empirically adjusting the distance between points (called \textit{radius} in the original paper) in order to generate meshes with the desired number of vertices.
Points are then connected in a Delaunay triangulation using the well known \textit{Triangle library} \cite{Shewchuk:2005}, with the default parameters configuration.\\
In $\Dtriangle$, $A_n$ and $e_n$ are almost constant, as no constraints are imposed to the triangulation process.

%%%%%%%%%%%%%%%%%%%%%%%%%%%%%%%%%%%%%
\subsection{Hybrid datasets}
\label{app:dataset:hybrid}
The construction of hybrid datasets is characterized by the insertion in $\Omega$ of one or more polygonal elements, and by a tessellation algorithm.
Each hybrid dataset is built around (and named after)  
an \textit{initial polygon} $\P=\P(t_n)$ depending on a \textit{deformation parameter} $t_n\in[0,1)$, which is used to deform $\P$.
This parameter directly depends on the mesh number (i.e. $t_n\to\ 1$ as $n\to N$), and it can be adjusted to improve or worsen the quality of the polygon $\P$ (the higher, the worse).\\
At refinement step $n$, mesh $\Omega_n$ is created by inserting a number of identical copies of the deformed polygon $\P(t_n)$ (opportunely resized) in the domain $\Omega$, and tessellating the rest of $\Omega$ using the Triangle library.
This procedure is detailed in Algorithm~\ref{alg:hybrid}.\\
Note that, a whole family of other datasets may be generated by simply defining a new initial polygon. More examples can be found in \cite{attene2019benchmark}.

\begin{figure}[htbp]
\centering
	\includegraphics[width=.6\linewidth]{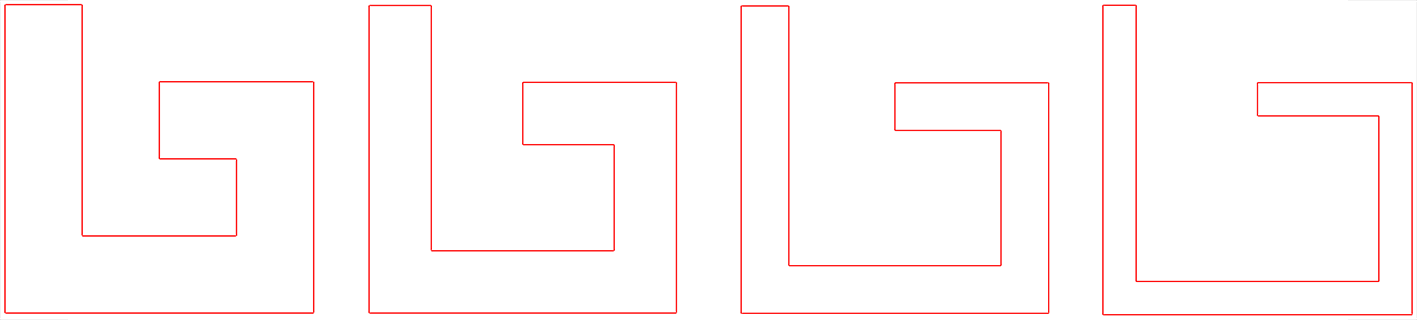} %% {Figures/basemesh-maze.png}
\caption{Initial polygons $\P(t_0)$, $\P(t_2)$, $\P(t_4)$ and $\P(t_6)$ from dataset $\Dmaze$.}
\label{fig:maze-single}
\end{figure}

\medskip
The initial polygon $\P(t_n)$ for dataset $\Dmaze$ is the 10-sided element shown in Fig.~\ref{fig:maze-single}, with vertices
\begin{align*}
    &(0,1),\ (0,0),\ (1,0),\ 
    (1,0.75),\ (0.5,0.75),\\
    &\left(0.5,0.5+\frac{t_n}{4}\right),\ \left(0.75+\frac{t_n}{4},0.5+\frac{t_n}{4}\right),\\
    &\left(0.75+\frac{t_n}{4},0.25-\frac{t_n}{4}\right),\ 
    \left(0.25-\frac{t_n}{4},0.25-\frac{t_n}{4}\right),\ 
    \left(0.25-\frac{t_n}{4},1\right).
\end{align*}
As $t_n\to 1$, the length of the shortest edge (the one with vertices $(0,1)$ and $(0.25-t_n/4,1)$) goes to zero, and so does the area of $\P(t_n)$.

\begin{figure}[htbp]
\centering
	\includegraphics[width=.6\linewidth]{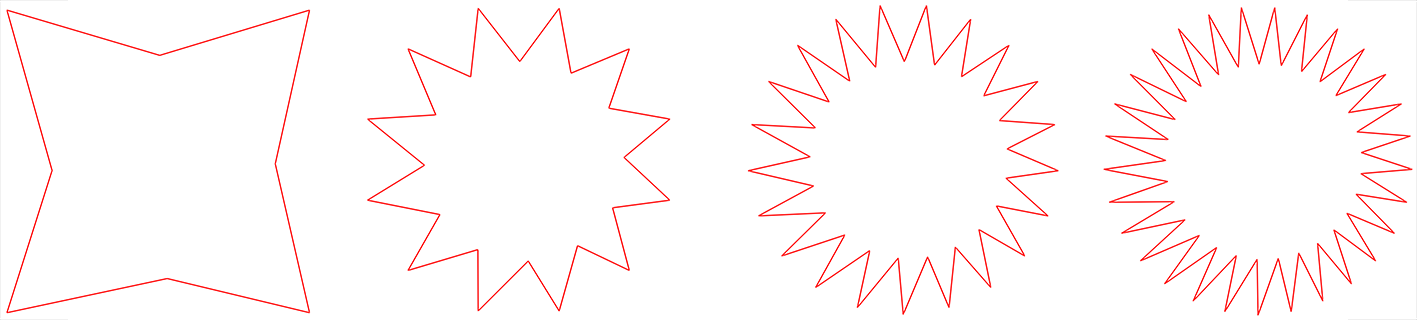} %% {Figures/basemesh-star.png}
\caption{Initial polygons $\P(t_0)$, $\P(t_2)$, $\P(t_4)$ and $\P(t_6)$ from dataset $\Dstar$.}
\label{fig:star-single}
\end{figure}

\medskip
For building the initial polygon $\P(t_n)$ of dataset $\Dstar$ (Fig.~\ref{fig:star-single}), we first build a $\Bar{i}$-sided regular polygon, with $\Bar{i}=8(1+\lfloor 10t_n \rfloor)$ and vertices
\begin{align*}
\begin{cases}
  v_0=(1,0),&\\
  v_i=\sigma(v_{i-1}), &\mbox{for $i=1, \ldots, \Bar{i}$},   
\end{cases}
\end{align*}
being $\sigma(v)$ the rotation centered at $(0,0)$ of vertex $v$ by an angle of $2\pi/\Bar{i}$.
Then we project every odd-indexed vertex towards the barycenter of $\P(t_n)$:
\begin{align*}
    v_{2j+1}'=s \ v_{2j+1}, \quad \mbox{for $j=0,\ldots,\frac{\Bar{i}-1}{2}$,}
\end{align*}
where the projection factor $s\in(0,1)$ is gradually decreased until the angles at the even-indexed vertices become smaller than $(1-t_n)\pi/3$.\\
As $t_n\to 1$ we have an increasing number of edges (from 8 to almost 90), the minimum angle and the area decrease to zero while the length of every edge increases.

\medskip
Once we defined the initial polygon $\P(t_n)$, we can build the corresponding dataset through Algorithm \ref{alg:hybrid}.
We have some initial parameters, which are set a priori and remain untouched: the number of meshes in the dataset $N$, the area of the initial polygon at the first step $d_0$ and the deformation range $T=[t_{\min}, t_{\max}]$.
In this work we set $N=10$, $d_0=0.03$, which corresponds to $3\%$ of the domain, and $T=[0,0.95]$.\\
Then we have three main parameters, $e_n\in \mathbb N$, $t_n\in T$ and $d_n\in (0,d_0)$, which respectively regulate the number of initial polygons inserted, the deformation of these polygons and their area.
In particular, $e_n$ increases inversely to $d_n$ ($\Omega_{n+1}$ has twice as polygons as $\Omega_n$, with halved areas), so that the percentage of the domain covered by polygons (not triangles) is preserved all across the dataset.
Due to the complicated shapes of some initial polygons, it may be hard to ask for exactly $|\P(t_n)|=d_n$, therefore we only impose $|\P(t_n)|\le d_n$.\\
Several options are possible for setting $e_n, t_n$ and $d_n$, and the speeds at which these quantities vary, strongly affect the geometrical qualities of the meshes in the dataset.
In our datasets, $e_n$ increases exponentially, $t_n$ increases linearly inside $T$ and $d_n$ decreases exponentially.
The exponential increase of the number of initial polygons inserted in the domain may lead to intersections between them, or with the domain boundaries.
To avoid this phenomenon, we inserted a \textit{while} loop in Algorithm~\ref{alg:hybrid} which decreases $d_n$ until no intersections occur: this ensures stability to the algorithm, but in practice it activates only for very dense meshes and it typically runs only few iterations.\\
Last, when all polygons have been inserted in $\Omega$, the Triangle algorithm is used to generate a Delaunay triangulation.
The already inserted polygons are considered as holes in the domain, and we set no limitations on the number of Steiner points that may appear in the triangulation process. We adopt the following parameters configuration, cf. \cite{Shewchuk:2005}:
\begin{itemize}
    \item \textbf{q:} no angles smaller than 20 degrees;
    \item \textbf{c:} enclose the convex hull with segments;
    \item \textbf{l:} use only vertical cuts in the divide-and-conquer algorithm (this switch is primarily of theoretical interest);
    \item \textbf{a:} maximum triangle area constraint, set equal to $d_n$.
\end{itemize}
Due to the freedom left to the Triangle algorithm, it is not possible to estimate $A_n$ and $e_n$ precisely for hybrid datasets; hence, the relative values reported in Table~\ref{table:ratios} have been measured a posteriori.
\begin{algorithm}[htbp]
\caption{hybrid datasets generation}
\begin{algorithmic}[1]
    \State define the initial polygon
    \State set the initial parameters $N$, $d_0$, and $T$
	\For {$n = 0, \ldots, N$}
		\State set the main parameters: $e_n=2^n$, $t_n=n\frac{t_{\max}-t_{\min}}{N}$, $d_n=d_0/2^n$
		\State use Poisson Disk Sampling with radius $r = 1/\sqrt{2e_n}$ to find a set of $e_n$ points $\{\bm{p_i}\}_{i=1,\ldots,e_n}$ in $\Omega$
		\State generate polygon $\P(t_n)$ with $|\P(t_n)|\le d_n$
		\State insert a copy of $\P(t_n)$ centered around every $\bm{p_i}$
		\While{polygon $\P(t_n)$ intersects with other polygons or with the boundary of $\Omega$}
		\State $d_n \leftarrow d_n-\epsilon$
	    \State generate a polygon $\P(t_n)$ with $|\P(t_n)|\le d_n$
		\State insert a copy of $\P(t_n)$ centered around every $\bm{p_i}$
		\EndWhile
		\State use Triangle to generate the Delaunay triangulation $\Omega_n$ of $\Omega$, considering polygons $\P(t_n)$ as holes
		\State add $\Omega_n$ to $\calD$
	\EndFor
\end{algorithmic} 
\label{alg:hybrid}
\end{algorithm}

%%%%%%%%%%%%%%%%%%%%%%%%%%%%%%%%%%%%%
\subsection{Mirroring datasets}
\label{app:dataset:mirroring}
The construction of $\Djenga$, $\Dslices$ and $\Dulike$, at every step $n\ge1$, consists in a first algorithm for iteratively generating a \textit{base mesh} $\widehat{\Omega}_n$ from the previous base mesh $\widehat{\Omega}_{n-1}$, followed by a \textit{mirroring technique} which returns the computational mesh $\Omega_n$.
The base mesh generation algorithm is different for each dataset (Algorithms \ref{alg:jenga}, \ref{alg:slices} and \ref{alg:ulike}), while the mirroring algorithm (Algorithm \ref{alg:mirror}) is common to all three datasets.
Algorithms \ref{alg:jenga}, \ref{alg:slices} and \ref{alg:ulike} depend on two initial parameters: $N$ indicates the number of meshes in the dataset and $N_{el}$ indicates the number of elements to insert at each step.
For mirroring datasets we set $N_{el}=1$, while for multiple mirroring datasets (described in the next section) we set $N_{el}=4$.

\begin{figure}[htbp]
\centering
	\includegraphics[width=.7\linewidth]{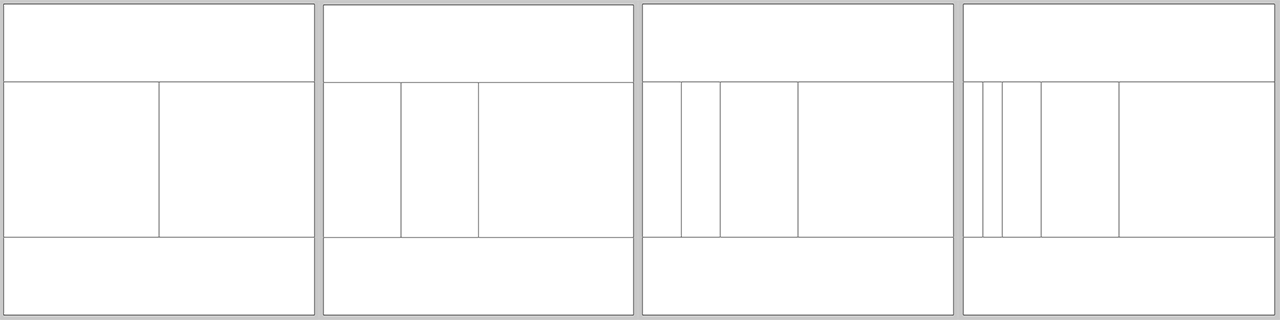} %% {Figures/basemesh-jenga.png}
\caption{Non-mirrored base meshes $\widehat{\Omega}_0, \widehat{\Omega}_1, \widehat{\Omega}_2, \widehat{\Omega}_3$ from datasets $\Djenga$.}
\label{fig:jenga-single}
\end{figure}

\medskip
In the $\Djenga$ base mesh shown in Fig.~\ref{fig:jenga-single} we have a \textit{top bar}, a \textit{bottom bar} and a \textit{right square} which are fixed independently of $n$, and $n+1$ \textit{rectangles} in the left part of the domain.
At each refinement step $n\ge 1$, a new rectangular element is created by splitting in two equal parts the leftmost rectangular element in the previous base mesh, and consequently updating the top and the bottom bars with new vertices and edges.
Therefore, all elements in $\widehat{\Omega}_n$, except for the top and the bottom bars, are rectangles with height equal to $1/2$ and basis ranging from $1/2$ to $1/2^{n+1}$.
Once that the base mesh $\widehat{\Omega}_n$ is generated, the mirroring algorithm is recursively applied for $n$ times to generate the computational mesh $\Omega_n$, as described in Algorithm~\ref{alg:jenga}.\\
When computing $A_n$ and $e_n$, we can restrict our calculations to the base mesh, because these ratios are not affected by the mirroring algorithm.
In particular, the longest edge in the base mesh is the upper edge of the top bar, which is never split, while the shortest edge is the basis of the leftmost rectangle, which halves at each step: this causes $e_n\sim 2^n$.
The top bar is also the element with the greatest area (together with the bottom bar and the right square), which is constantly equal to 1/4, while the leftmost rectangle has area $1/2*1/2^{n+1}=1/2^{n+2}$, therefore $A_n\sim 2^n$.\\

\begin{algorithm}[htbp]
\caption{$\Djenga$ dataset generation}
\begin{algorithmic}[1]
    \State set the number of meshes $N$ and the number of elements $N_{el}$
    \For{$n=0, \ldots, N$}
        \State \textit{top bar} = $\{(0,0.75), (1, 0.75), (1, 1), (0, 1)\}$
        \State \textit{bottom bar} = $\{(0,0), (1, 0), (1, 0.25), (0, 0.25)\}$
        \State \textit{right square} = $\{(0.5,0.25), (1, 0.25), (1, 0.75), (0.5, 0.75)\}$
        \State vector $\boldsymbol{b}$ = sample $n*N_{el}$ equally spaced points inside interval $(0,0.5)$
        \For{$i=1, \ldots, \text{size}(\boldsymbol{b})$}
            \State \textit{rectangles}$[i]$ = $\{(\boldsymbol{b}[i-1], 0.25), (\boldsymbol{b}[i], 0.25), (\boldsymbol{b}[i], 0.75), (\boldsymbol{b}[i-1], 0.75)\}$
            \State insert vertex $(\boldsymbol{b}[i], 0.75)$ in \textit{top bar}
            \State insert vertex $(\boldsymbol{b}[i], 0.25)$ in \textit{bottom bar}
        \EndFor
        \State generate mesh $\widehat{\Omega}_n=\{\textit{top bar, bottom bar, right square, rectangles}\}$ 
        \For{$i=1, \ldots, n$}
            \State $\widehat{\Omega}_n$ = mirror mesh$(\widehat{\Omega}_n)$
        \EndFor
	\State add the newly generated $\Omega_n=\widehat{\Omega}_n$ to $\Djenga$
    \EndFor
\end{algorithmic} 
\label{alg:jenga}
\end{algorithm}

\begin{figure}[htbp]
\centering
	\includegraphics[width=.7\linewidth]{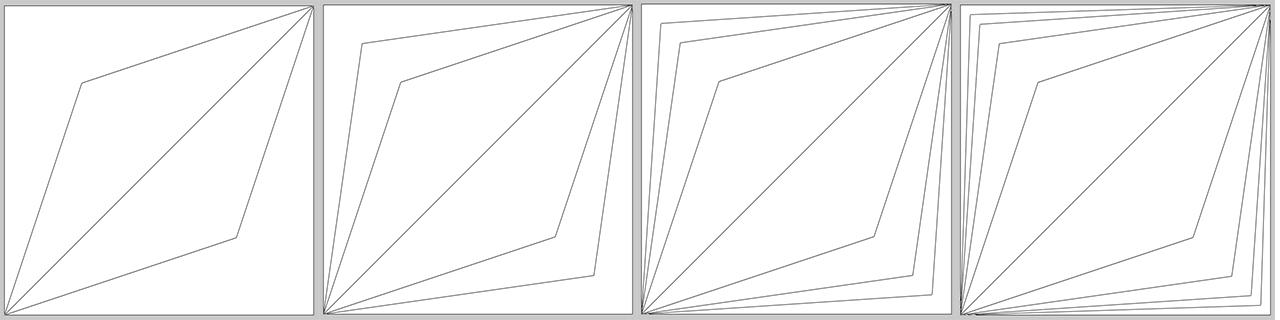} % {Figures/basemesh-slices.png}
\caption{Non-mirrored base meshes $\widehat{\Omega}_0, \widehat{\Omega}_1, \widehat{\Omega}_2, \widehat{\Omega}_3$ from datasets $\Dslices$.}
\label{fig:slices-single}
\end{figure}

\medskip
In the $\Dslices$ base meshes shown in Fig.~\ref{fig:slices-single}, at each step $n\ge 0$, we add the vertices with coordinates $(2^{-i}, 1-2^{-i})$ and $(1-2^{-i}, 2^{-i})$ for $i=1,\ldots ,n+2$, and we connect them to the vertices $(0,0)$ and $(1,1)$.
As a result, at each iteration we create a couple of new polygons, called \textit{upper slice} and \textit{lower slice}, symmetrical with respect to the diagonal, and we add them to the base mesh.\\
The area of the two inner triangles (the biggest polygons in the base mesh) is constantly equal to $1/4$.
For evaluating the area of the two most external polygons, we consider them as the union of the two identical triangles obtained by splitting the polygons along the diagonal (the one connecting the vertices with coordinates $(0,1)$ and $(1,0)$).
Then the smallest area in the base mesh is the sum of the areas of two equal triangles with basis $\sqrt{2}/2$ and height $2^{-n}/\sqrt{2}$, and simple calculations lead to $A_n \sim 2^{n}$.
Last, we notice that all the edges in the base mesh have lengths between 1 and $\sqrt{2}$, because no edge is ever split, hence $e_n \sim c$.
\begin{algorithm}[htbp]
\caption{$\Dslices$ dataset generation} 
\begin{algorithmic}[1]
    \State set the number of meshes $N$ and the number of elements $N_{el}$
    \For{$n=0, \ldots, N$}
        \State vector $\boldsymbol{b}$ = $[2^{-1}, 2^{-2}, \ldots, 2^{-{n*N_{el}}}]$
        \For{$i=1, \ldots, \text{size}(\boldsymbol{b})$}
            \State \textit{upper slices}$[i]$ = $\{(0, 0), (\boldsymbol{b}[i], 1-\boldsymbol{b}[i]), (1, 1), (\boldsymbol{b}[i+1], 1-\boldsymbol{b}[i+1])\}$
            \State \textit{lower slices}$[i]$ = $\{(0, 0), (1-\boldsymbol{b}[i], \boldsymbol{b}[i]), (1, 1), (1-\boldsymbol{b}[i+1], \boldsymbol{b}[i+1])\}$
        \EndFor
        \State generate mesh $\widehat{\Omega}_n=\{\textit{upper slices, lower slices}\}$ 
        \For{$i=1, \ldots, n$}
            \State $\widehat{\Omega}_n$ = mirror mesh$(\widehat{\Omega}_n)$
        \EndFor
	\State add the newly generated $\Omega_n=\widehat{\Omega}_n$ to $\Dslices$
    \EndFor
\end{algorithmic} 
\label{alg:slices}
\end{algorithm}

\begin{figure}[htbp]
\centering
	\includegraphics[width=.7\linewidth]{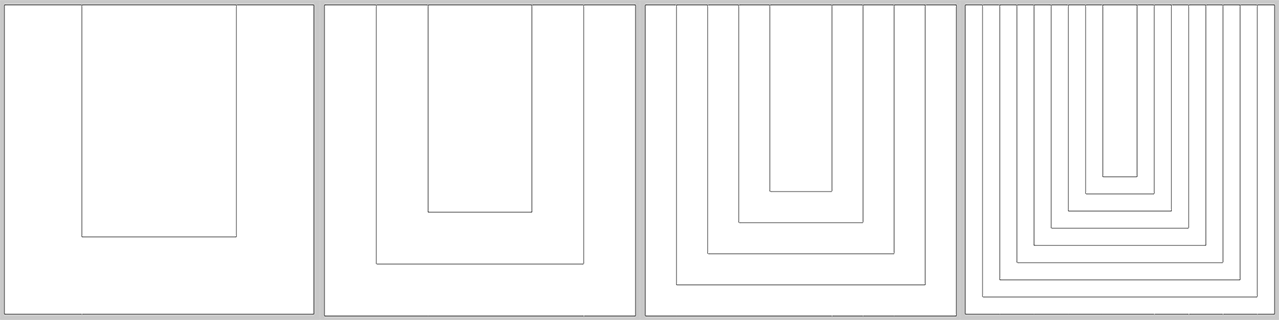} %% {Figures/basemesh-ulike.png}
\caption{Non-mirrored base meshes $\widehat{\Omega}_0, \widehat{\Omega}_1, \widehat{\Omega}_2, \widehat{\Omega}_3$ from datasets $\Dulike$.}
\label{fig:ulike-single}
\end{figure}

\medskip
In the $\Dulike$ base meshes shown in Fig.~\ref{fig:ulike-single}, at each step $n\ge0$ we insert $2^n$ $U$-shaped continuous polylines inside the domain.
We have an \textit{internal rectangle} and a sequence of concentric equispaced \textit{U-like polygons} culminating with the \textit{external U}.
This last element is not different from the other $U$-like polygons, but is created separately, because we need to split its lower edge in order to match the base mesh that will appear below it during the mirroring algorithm.\\
In every base mesh, the shortest edge is the one corresponding to the width of each $U$-like polygon, which measures $2^{-(n+1)}$, and the longest edges are the left and right boundaries of the domain.
This causes $e_n \sim 2^{n}$.
Said $e$ the shortest edge, the smallest area is the one of the internal rectangle, equal to $2e(1/2+e)$, and the biggest area is the one relative to the external $U$, equal to $3e-2e^2$.
We have
\begin{align*}
    A_n = \frac{3-2e}{1+2e} = \frac{3-2(2^{-(n+1)})}{1+2(2^{-(n+1)})} = \frac{3-2^{-n}}{1+2^{-n}}\sim c.
\end{align*}

\begin{algorithm}[htbp]
\caption{$\Dulike$ dataset generation}
\begin{algorithmic}[1]
    \State set the number of meshes $N$ and the number of elements $N_{el}$
    \For{$n=0, \ldots, N$}
        \State vector $\boldsymbol{b}$ = sample $2^{n*N_{el}}$ equally spaced points inside interval $(0,0.5)$
        \For{$i=1, \ldots, \text{size}(\boldsymbol{b})$}
            \State \textit{U-like polygons}$[i]$ = $\{
            (\boldsymbol{b}[i], 1), 
            (\boldsymbol{b}[i], \boldsymbol{b}[i]), 
            (1-\boldsymbol{b}[i], \boldsymbol{b}[i]),
            (1-\boldsymbol{b}[i], 1), 
            (1-\boldsymbol{b}[i+1], 1), 
            (1-\boldsymbol{b}[i+1], \boldsymbol{b}[i+1]),
            (\boldsymbol{b}[i+1], \boldsymbol{b}[i+1]), 
            (\boldsymbol{b}[i+1], 1)\}$
        \EndFor
        \State \textit{internal rectangle} = $\{(\boldsymbol{b}[end], 1), (\boldsymbol{b}[end], \boldsymbol{b}[end]), (1-\boldsymbol{b}[end], \boldsymbol{b}[end]), (1-\boldsymbol{b}[end], 1)\}$
        \State \textit{external U} = $\{(0,1), (0,0), (1,0), (1,1), (1-\boldsymbol{b}[0], 1), (1-\boldsymbol{b}[0], \boldsymbol{b}[0]), (\boldsymbol{b}[0], \boldsymbol{b}[0]), (\boldsymbol{b}[0], 1)\}$
        \For{$b \in \boldsymbol{b}$}
            \State insert vertices $(b,0)$ and$(1-b,0)$ in \textit{external U}
        \EndFor
        \State generate mesh $\widehat{\Omega}_n=\{\textit{external U, U-like polygons, internal rectangle}\}$ 
        \For{$i=1, \ldots, n$}
            \State $\widehat{\Omega}_n$ = mirror mesh$(\widehat{\Omega}_n)$
        \EndFor
	\State add the newly generated $\Omega_n=\widehat{\Omega}_n$ to $\Dulike$
    \EndFor
\end{algorithmic} 
\label{alg:ulike}
\end{algorithm}

%%%%%%%%%%%%%%%%%%%%%%%%%%%%%%%%%%%%%
\subsection{Multiple mirroring datasets}
\label{app:dataset:multiple-mirroring}
\begin{figure}[htbp]
\centering
	\includegraphics[width=.55\linewidth]{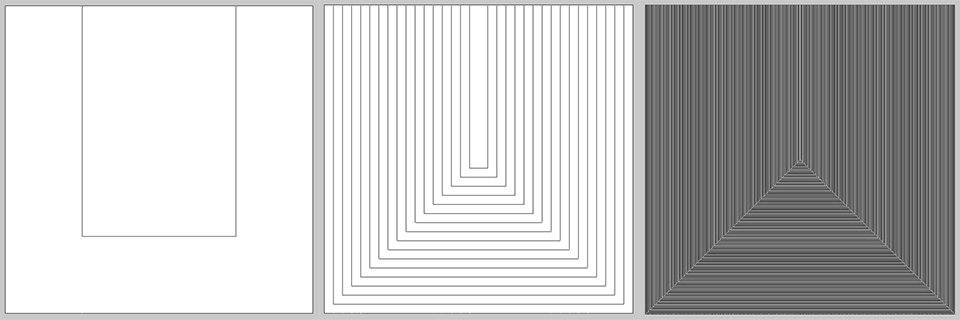} %% {Figures/basemesh-ulike4.png}
\caption{Non-mirrored base meshes $\widehat{\Omega}_0, \widehat{\Omega}_1$ and $\widehat{\Omega}_2$ from datasets $\DulikeM$.}
\label{fig:ulike-single-x4}
\end{figure}

Multiple mirroring datasets are built with the exactly same algorithms of the mirroring datasets, changing the parameter $N_{el}$.
This parameter regulates the number of elements generated in each base mesh of the dataset.
In particular, datasets $\DjengaM$, $\DslicesM$ and $\DulikeM$ are defined setting $N_{el}=4$.
An example of a multiple mirroring dataset with $N_{el}=4$ is shown in Fig.~\ref{fig:ulike-single-x4}, where the first three base meshes of $\DulikeM$ are presented.\\
The $N_{el}$ value influences ratios $A_n$ and $e_n$: if $A_n, e_n\sim 2^n$ for $N_{el}=1$, these quantities become asymptotic to $2^{4n}$ when $N_{el}=4$, except for the cases in which the ratios were constant (see Table~\ref{table:ratios}).

%%%%%%%%%%%%%%%%%%%%%%%%%%%%%%%%%%%
\subsection{The mirroring algorithm}
\label{app:dataset:mirroring-algorithm}
The \textit{mirroring} algorithm (Algorithm~\ref{alg:mirror}) generates four adjacent copies of any polygonal mesh $\calM$ defined over the domain $\Omega=[0,1]^2$.
In CinoLib \cite{livesu2019cinolib}, a polygonal mesh can be defined by a vector \textit{verts} containing all its vertices and a vector \textit{polys} containing all its polygons.
The result of the algorithm is therefore a polygonal mesh $\calM'$, generated by some vectors \textit{new-verts} and \textit{new-polys}, containing four times the number of vertices and polygons of $\calM$.
When iterated a sufficient number of times, this construction allows us to obtain a number of vertices and degrees of freedom in each mesh of the mirroring datasets that is comparable to that of the meshes at the same refinement level in hybrid datasets.\\
Vector \textit{new-verts} contains all vertices $v\in$ \textit{verts} copied four times and translated by vectors $(0,0)$, $(1,0)$, $(1,1)$ and $(0,1)$ respectively.
The coordinates of all vertices in \textit{new-vertices} are divided by 2, so that all new points lie in the same domain as the input mesh.
Vector \textit{new-polys} is simply vector \textit{polys} repeated four times.
A final cleaning step is required to remove duplicated vertices and edges that may arise in the mirroring process, for example if the initial mesh $\calM$ has vertices along its boundary.
\begin{algorithm}[htbp]
\caption{mesh mirroring} 
\begin{algorithmic}[1]
    \State input: base mesh $\calM$
    \State \textit{verts} = verts($\calM$), \ \textit{polys} = polys($\calM$)
    \State \textit{new-verts} = \textit{verts}
    \For{vertex $v \in$ \textit{verts}}
    \State insert vertex $v+(1,0)$ in \textit{new-verts}
    \EndFor
    \For{vertex $v \in$ \textit{verts}}
    \State insert vertex $v+(1,1)$ in \textit{new-verts}
    \EndFor
    \For{vertex $v \in$ \textit{verts}}
    \State insert vertex $v+(0,1)$ in \textit{new-verts}
    \EndFor
    \For{vertex $v \in$ \textit{new-verts}}
    \State $v \leftarrow v/2$
    \EndFor
    \State \textit{new-polys} = [\textit{polys}, \textit{polys}, \textit{polys}, \textit{polys}]
    \State $\calM'$ = mesh$\{\textit{new-verts}, \textit{new-polys}\}$
    \State remove duplicated vertices and edges from $\calM'$
\end{algorithmic} 
\label{alg:mirror}
\end{algorithm}

\section*{Acknowledgements}
This paper has been realised in the framework
of ERC Project CHANGE, which has received funding from the European Research Council (ERC) under the European Union’s Horizon 2020 research and innovation program (grant agreement no.~694515).
We are very grateful to Dr. L. Mascotto, University of Vienna, for useful suggestions.

%%\bibliographystyle{plain}
%%\bibliography{vem}

\end{document}